\documentclass[hidelinks]{article}
\newcommand{\documentdate}{9 VII 2025}

\usepackage{a4wide,latexsym,graphicx,amsmath,amsfonts}
\usepackage{varioref}
\usepackage{xcolor}
\usepackage{ulem}
\usepackage{multirow}
\usepackage{hyperref}
\title{Recursive Bound-Constrained AdaGrad with Applications\\ 
to Multilevel and Domain Decomposition Minimization}

\author{Serge Gratton\footnotemark[1],
        Alena Kopani\v{c}\'akov\'a\footnotemark[2],
        Philippe L. Toint\footnotemark[3]}

\newcommand{\beqn}[1]{\begin{equation}\label{#1}}
\newcommand{\eeqn}{\end{equation}}
\newcommand{\req}[1]{(\ref{#1})}
\newcommand{\ms}{\;\;\;\;}
\newcommand{\tim}[1]{\;\; \mbox{#1} \;\;}

\newcounter{algo}[section]
\renewcommand{\thealgo}{\thesection.\arabic{algo}}
\newcommand{\llem}[2]{\vspace{\baselineskip} 
\noindent\framebox[\textwidth]{\parbox{0.95\textwidth}{
\begin{lemma} \label{#1} \rm #2 \end{lemma} } } \vspace{\baselineskip} }
\newcommand{\algo}[3]{\refstepcounter{algo}
\begin{center}\begin{figure}[htbp]
\framebox[1.05\textwidth]{
\parbox{\textwidth} {\vspace{\topsep}
{\bf Algorithm \thealgo : #2}\label{#1}\\
\vspace*{-\topsep} \mbox{ }\\
{#3} \vspace{\topsep} }}
\end{figure}\end{center}}
\newcommand{\bpr}{{\bf Proof.} \hspace{1.5mm}}
\newcommand{\epr}{\hfill $\Box$ \vspace*{1em}}
\newcommand{\lthm}[2]{\vspace{\baselineskip} 
\noindent\framebox[\textwidth]{\parbox{0.95\textwidth}{
\begin{theorem} \label{#1} \rm #2 \end{theorem} } } \vspace{\baselineskip} }
\newcommand{\lcor}[2]{\vspace{\baselineskip} 
\noindent\framebox[\textwidth]{\parbox{0.95\textwidth}{
\begin{corollary} \label{#1} \rm #2 \end{corollary} } } \vspace{\baselineskip}
}
\newcommand{\ii}[1]{\{ 1, \ldots, #1 \}}
\newcommand{\iiz}[1]{\{ 0, \ldots, #1 \}}
\newcommand{\iibe}[2]{\{ #1, \ldots, #2 \}}
\DeclareMathAlphabet{\pazocal}{OMS}{zplm}{m}{n}
 
\newcommand{\calB}{{\pazocal{B}}} 
\newcommand{\calC}{{\pazocal{C}}} 
\newcommand{\calD}{{\pazocal{D}}} 
\newcommand{\calE}{{\pazocal{E}}} 
\newcommand{\calF}{{\pazocal{F}}} 
\newcommand{\calI}{{\pazocal{I}}}

\newcommand{\calO}{{\pazocal{O}}} 
\newcommand{\calR}{{\pazocal{R}}}

\renewcommand{\Re}{\hbox{I\hskip -2pt R}}

\newcommand{\bigfrac}[2]{\frac{\displaystyle #1}{\displaystyle #2}}

\newcommand{\sfrac}[2]{{\scriptstyle \frac{#1}{#2}}}
\newcommand{\half}{\sfrac{1}{2}}
\newcommand{\eqdef}{\stackrel{\rm def}{=}}
\newcommand{\bigsum}{\displaystyle \sum}
\newcommand{\bigsqrt}{\displaystyle \sqrt}
\newcommand{\kap}[1]{\kappa_{\mbox{\tiny #1}}}

\newcommand{\flow}{f_{\rm low}}
\newcommand{\al}[1]{{\footnotesize{\sf #1}}}

\newtheorem{theorem}{Theorem}[section]
\newtheorem{lemma}[theorem]{Lemma}
\newtheorem{corollary}[theorem]{Corollary}

\newcommand{\adagrad}{\al{ADAGB2}}
\newcommand{\algnameSL}{\al{ADAGB2}}
\newcommand{\algname}{\al{ML-ADAGB2}}
\newcommand{\algnamed}{\al{DD-ADAGB2}}
\newcommand{\algnamehyb}{\al{ML-DD-ADAGB2}}
\newcommand{\proof}[1]{
\begin{list}{}{
\setlength{\topsep}{0.0pt}
\setlength{\partopsep}{0.0pt}
\setlength{\leftmargin}{0.025\textwidth}
\setlength{\rightmargin}{0.5\leftmargin}
\setlength{\labelwidth}{0.5\leftmargin}
\setlength{\labelsep}{0.25\leftmargin}}
\item \bpr #1 \epr \noindent
\end{list}}
\DeclareMathOperator*{\average}{average}

\newcommand{\EA}[1]{\mathbb{E}^\calE\!\left[#1 \right]}

\newcommand{\EcondA}[2]{\mathbb{E}_{#1}^\calE\!\left[#2 \right]}
\newcommand{\Prob}[1]{\mathbb{P}\!\left[#1 \right]}

\definecolor{myblack}{RGB}{53, 53, 53}
\definecolor{myblue}{RGB}{40, 75, 99}
\definecolor{myred}{RGB}{192, 50, 33}
\definecolor{myyellow}{RGB}{255, 166, 48}
\definecolor{mywhite}{RGB}{240, 237, 238}
\definecolor{mygreen}{RGB}{0, 102, 0}

\definecolor{green1}{RGB}{9, 82, 86}
\definecolor{green2}{RGB}{8, 127, 140}
\definecolor{green3}{RGB}{6, 167, 125}
\definecolor{green4}{RGB}{79, 109, 122}
\definecolor{green5}{RGB}{192, 214, 223}
\definecolor{violet}{RGB}{26,69,131}

\definecolor{checkgreen}{rgb}{0,0.6,0}
\definecolor{phase1}{rgb}{0.008,0.655,1.000}
\definecolor{phase2}{rgb}{0.016,0.75,0.700}
\definecolor{phase3}{rgb}{0.929,0.35,0.700}
\definecolor{icsyellow}{cmyk}{0.00,0.11,0.53,0.00}

\usepackage{tikz}
\usepackage{pgfplotstable}
\usepackage{pgfplots}
\usepackage[outline]{contour}
\usepackage{calrsfs}
\usepackage{graphicx}
\usepgfplotslibrary{groupplots}
\pgfplotsset{width=7cm,compat=1.3}
\usetikzlibrary{matrix}

\usetikzlibrary{external}
\date{\documentdate}

\begin{document}

\maketitle

\renewcommand{\thefootnote}{\fnsymbol{footnote}}
\footnotetext[1]{Toulouse-INP, IRIT, ANITI, Toulouse,
  France. Email: serge.gratton@toulouse-inp.fr.}
\footnotetext[2]{Toulouse-INP, IRIT, ANITI, Toulouse,
  France. Email: alena.kopanicakova@toulouse-inp.fr.}
\footnotetext[3]{Namur Center for Complex Systems (naXys),
  University of Namur, Namur, Belgium. Email: philippe.toint@unamur.be.}

\renewcommand{\thefootnote}{\arabic{footnote}}

\begin{abstract}
Two OFFO (Objective-Function Free Optimization) noise tolerant algorithms are presented that handle bound constraints, inexact gradients and use second-order information when available.  The first is a multi-level method exploiting a hierarchical description of the problem and the second is a domain-decomposition method covering the standard addditive Schwarz decompositions.  Both are generalizations of the first-order AdaGrad algorithm for unconstrained optimization. Because these algorithms share a common theoretical framework, a single convergence/complexity theory is provided which covers them both. Its main result is that, with high probability, both methods need at most $\calO(\epsilon^{-2})$ iterations and noisy gradient evaluations to compute an $\epsilon$-approximate first-order critical point of the bound-constrained problem.  Extensive numerical experiments are discussed on applications ranging from PDE-based problems to deep neural network training, illustrating their remarkable computational efficiency.
\end{abstract}


\section{Introduction}

Consider the problem
\beqn{problem}
\min_{x\in \calF} f_r(x),
\tim{ where }
\calF = \big\{x\in\Re^{n_r} \mid l_{r,i} \leq x_i \leq u_{r,i} \tim{for} i \in \ii{n_r} \big\}
\eeqn
and where $f_r$ is a smooth possibly nonconvex 
and noisy function from $\Re^n$ into $\Re$. 
Problems of this type frequently arise in diverse 
scientific applications, e.g.,~contact or fracture 
mechanics~\cite{kopanivcakova2020recursive, krause2009nonsmooth}. 
Another prominent application area involves inverse 
problems and PDE-constrained optimization (under 
uncertainty)~\cite{ciaramella2024multigrid}, where 
the reduced functional is often subject to bounds 
ensuring the physical feasibility.
Moreover, in recent years, bound-constrained 
minimization problems have also emerged in machine 
learning due to the incorporation of prior knowledge 
into the design of machine learning models.
For instance, features or 
weights~\cite{leimkuhler2021better} are often 
restricted to specific ranges to prevent 
overfitting. Moreover, in scientific machine 
learning, bounds are frequently applied to reflect 
real-world limits and maintain physical 
interpretability~\cite{lu2021physics}.
Although these application areas are fairly diverse, 
the arising minimization problems face a common 
challenge: solving the underlying minimization 
problems is computationally demanding due to the 
typically large dimensionality $n \in \mathbb{N}$, 
severe ill-conditioning and the presence of both noise and bounds.

Noise in the objective function may be caused by a number of factors (such as data sampling or variable precision computations, for instance), but we subsume all these cases here
by considering that noise affects the objective function randomly. Various strategies have been proposed to handle this situation in the standard  linesearch, trust-region or adaptive regularization settings (see \cite{Cart91,CartSche17,BellKrejKrkl18,BlanCartMeniSche19,GratToin20,BellGuriMoriToin21b,BeraCaoSche21,BellGuriMoriToin22}, to cite only a few), but one of the most successful is the Objective-Function-Free-Optimization (OFFO) technique.  Taking into account that optimization of noisy functions require more accuracy on the objective function value than on its derivatives, OFFO algorithms bypass this difficulty by never computing such values, yielding substantially improved reliability in noisy situations. Popularized in (unconstrained) machine learning, where noise is typically caused by sampling, OFFO algorithms such as stochastic gradient descent (SGD) \cite{RobbMunr51}, AdaGrad~\cite{DuchHazaSing11} or Adam \cite{KingBa15} (among many others) have had a significant impact on numerical optimization.

In this paper, we propose \algname\ and \algnamed, two new OFFO variants of \adagrad\  in order to 
alleviate the computational challenge mentioned above. Both variants allow the use of second-order information, should it be available, bounds on the variables and stochastic noise on the gradients, while exploiting structure frequently present in problems of the form \req{problem}.

A first example of such structure occurs 
when~\eqref{problem} arises from the 
discretization of an infinite-dimensional problem, in 
which case cheap surrogates of the objective 
function can be constructed by 
exploiting discretizations with lower resolutions. 
Multilevel methods are well known to take advantage 
of such structure. They have been originally 
proposed to solve discretized elliptic partial 
differential equations (PDEs)~\cite{BrigHensMcCo00}, 
for which they have been show to exhibit optimal 
complexity and convergence rate independent of the 
problem size. Since then, the multilevel approaches 
have been successfully extended to solve non-convex 
minimization problems, using multilevel line-search 
(MG-OPT)~\cite{Nash00}, multilevel trust-region 
(RMTR)~\cite{GratSartToin08} or multilevel higher-
order regularization strategies 
(M-ARC)~\cite{calandra2021high}. 

Extending multilevel methods to handle bound-constraints is challenging, as the lower levels are 
often unable to resolve fine-level constraints 
sufficiently well, especially when the constraints 
are oscillatory~\cite{kornhuber2001adaptive}.
Initial attempts to incorporate constraints into the 
multilevel framework involved solving linear 
complementarity problems; see, for 
instance,~\cite{mandel1984multilevel, brandt1983multigrid, hackbusch1983multi, gelman1990multilevel}.
These methods employed various constraint projection 
rules to construct lower-level bounds so 
that prolongated lower-level corrections would not 
violate the finest level bounds. 
However, these projection rules tend to be overly 
restrictive, resulting in multilevel methods that 
converge significantly slower than standard linear 
multigrid. To improve convergence speed, Kornhuber
proposed a truncated monotone multigrid 
method~\cite{kornhuber1994monotone} which employs 
a truncated basis approach and recovers the 
convergence rate of the unconstrained multigrid once 
the exact active-set is detected~\cite{hoppe1994adaptive, 
kornhuber1994monotone}.

In the context of nonlinear optimization, only a few 
multilevel algorithms for bound constrained problems 
exist. For instance, Vallejos proposed a gradient 
projection-based multilevel 
method~\cite{vallejos2010mgopt} for solving bilinear 
elliptic optimal control problems. 
In~\cite{kovcvara2016first}, two multilevel linesearch methods were introduced for convex 
optimization problems, incorporating the constraint 
projection rules from~\cite{hackbusch1983multi} and 
an active-set approach 
from~\cite{kornhuber1994monotone}.
For non-convex optimization problems, Youett et al.~proposed a filter trust-region 
algorithm~\cite{youett2019globally} that employs an 
active-set multigrid 
method~\cite{kornhuber1994monotone} to solve the 
resulting linearized problems. Additionally, Gratton et al.~introduced a variant of 
the RMTR method~\cite{GratMoufToinWebe08}, utilizing 
constraint projection rules from~\cite{gelman1990multilevel}. 
This approach was later enhances 
in~\cite{kopanivcakova2023multilevel} by employing 
the Kornhuber's truncated basis approach.
Our new \algname\ algorithm follows similar ideas, but in the noise tolerant OFFO context.

A second important example of problem structure, again typically arising in problems 
associated with PDEs, involves cases where the problem’s variables pertain to subdomains 
on which the underlying problem can be solved at a reasonable cost. 
This naturally leads to a class of methods known as domain-decomposition methods.
As with multilevel methods, domain-decomposition methods were originally developed for 
solving elliptic PDEs, giving rise to techniques such as the (restricted) additive Schwarz 
methods~\cite{dolean2015introduction,CaiSark99}. In the context of nonlinear problems, 
much of the focus has been on improving the convergence of Newton’s method. 
In particular, Cai et al.~proposed an additively preconditioned inexact Newton method, 
known as ASPIN~\cite{cai2002nonlinearly}. The core idea behind this approach is to 
rebalance the nonlinearities by solving restricted nonlinear systems on subdomains. 
Since then, numerous variants of additive preconditioners for Newton's method have 
appeared in the literature, e.g.,~RASPEN~\cite{dolean2016nonlinear} and 
SRASPEN~\cite{chaouqui2022linear}.

Within the framework of nonlinear optimization, the literature is comparatively sparse, 
but several inherently parallel globalization strategies have been developed. These 
include parallel variable and gradient 
distribution techniques~\cite{ferris1994parallel,mangasarian1995parallel}, additively preconditioned 
line search methods~\cite{park2022additive}, and domain-decomposition-based trust-region 
approaches~\cite{gross2009unifying, gross2009new, GratViceZhan19}.
More recently, several efforts have explored the use of nonlinear domain decomposition methods for 
machine learning applications, such as additively preconditioned 
L-BFGS~\cite{kopanivcakova2024enhancing} or trust-region method~\cite{trotti2025parallel}. 
While these approaches show promise for improving training efficiency and increasing model accuracy, to 
the best of our knowledge, none of the existing domain decomposition algorithms have yet been designed 
to be tolerant to (subsampling) noise, as required in practical applications.

Extending domain-decomposition methods to handle bound constraints is generally simpler than in the case 
of multilevel methods, since constraints may often be easily restricted to each subdomain.
To this aim, Badia et al.~developed the variants of additive Schwarz for linear variational 
inequality~\cite{badea2000additive,badea2004convergence,badea2012one}. 
In~\cite{kothari2022nonlinear}, a nonlinear restricted additive Schwarz preconditioner was introduced to 
precondition a Newton-SQP algorithm.
More recently, Park proposed the additive Schwarz framework for solving the convex optimization problems 
based on the (accelerated) gradient methods~\cite{park2020additive, park2021accelerated}. 
In the context of non-convex optimization, trust-region-based domain decomposition algorithms that 
enforce the trust-region constraint in the infinity norm naturally facilitate the incorporation of bound 
constraints; see, for example,~\cite{gross2009new}.

As noted in~\cite{GratViceZhan19,GratMercRiccToin23}, domain decomposition methods share strong 
conceptual similarities with multilevel approaches. 
Indeed, combining multilevel and domain decomposition methods is often essential in order to design 
highly parallel and algorithmically scalable algorithms. 
One of the key objectives of this work is to leverage this connection to develop a new domain 
decomposition method, \algnamed, which—like \algname—operates without requiring function values, 
supports bound constraints, and remains robust under stochastic gradient approximations.

We summarize our contributions as follows.
\begin{enumerate}
    \item We propose the OFFO multilevel \algname\ algorithm for problem exhibiting a hierarchical description, whose distinguishing features are its capacity to handle bound constraints and its good reliability in the face of noise (due to its OFFO nature and its tolerance to random noise in the gradients).  It is also capable of exploiting second-order information, should it be accessible at reasonable cost.
\item We then apply the concepts developed for 
the multilevel case to the context of domain-decomposition and propose the more 
specialized \algnamed\ algorithm with the same features.
\end{enumerate}
These contributions build upon the analysis of 
\cite{BellGratMoriToin25} for the single level stochastic OFFO case, 
the techniques presented in \cite{GratSartToin08,GratMoufToinWebe08} 
for multilevel deterministic non-OFFO and in \cite{GratKopaToin23} for unconstrained multilevel OFFO,
and on the discussion of the domain-decomposition framework in \cite[Section~5]{GratViceZhan19}.

After this introduction, Section~\ref{section:algo} describes hierarchical problems 
in more detail and proposes the \algname\ OFFO algorithm.  Its convergence and 
computational complexity is then analyzed in Section~\ref{section:convergence}, while 
Section~\ref{section:decomposition} describes its extension to the domain-decomposition 
context. Numerical experiments presented in Section~\ref{section:numerics} illustrate 
the efficiency, scope and versatility of the proposed methods. Some conclusions are 
finally discussed in Section~\ref{section:conclusion}.

\section{An OFFO multilevel algorithm for hierarchical problems}
\label{section:algo}

Given problem \req{problem}, we assume that there exists a 
hierarchy of spaces $\{\Re^{n_\ell}\}_{\ell=0}^r$ where (possibly local or simplified) version of the 
problem can be described\footnote{Note that we have not prescribed that $n_{\ell-1}\leq n_\ell$.}.
We refer to these spaces as \textit{levels}, from the top one (indexed by $\ell=r)$ 
to the bottom one (indexed by $\ell=0$).
Our algorithm is iterative and generates iterates $x_{\ell,k} \in \Re^{n_\ell}$, where 
the first subscripts denotes the level and the second the iteration within that level.  
At the top level, the task consists in minimizing $f_r$, the problem's objective 
function. For each feasible iterate $x_{r,k}$ it is possible to construct, at level 
$r-1$, a model $f_{r-1}$ of $f_r$, using local information such has $g_{r,k}$ 
(an approximation of $\nabla_x f_r(x_{r,k})$) and the \textit{restriction} 
operator $R_r$.  The bounds on the variables, if present, can also be transferred 
to level $r-1$ using the $R_r$, yielding a feasible set at level $r-1$, denoted 
$\calF_{r-1}$. The minimization may then be pursued at level $r-1$, descending on 
$f_{r-1}$ and generating iterates 
${x_{r-1,0}, x_{r-1,1}, \ldots  \ \text{in} \  \calF_{r-1}}$. Again, at each of 
these iterates, it is possible to construct, at level $r-2$, a model $f_{r-2}$ of 
$f_{r-1}$ and a feasible set $\calF_{r-2}$, and so on until level $l=0$ is reached, 
where the construction of a lower level approximation of $f_0$ is either impossible 
or too coarse.  Once an iteration at level $\ell$ is concluded (either because 
an approximate critical point is reached or because the maximum number of iterations 
is exceeded), the final iterate $x_{\ell,*}$ is \textit{prolongated} to level 
$\ell+1$ using the operator $P_{\ell+1}$.

As we have already mentioned, the iterates at level $\ell$ are denoted by $x_{\ell,k}$ 
for values of $k$ increasing from 0. Note that this is slightly ambiguous since iterates 
generated at level $\ell-1$ from different iterates at level $\ell$ share the same 
indices, but a fully explicit notation is too cumbersome, and no confusion arises in 
our development.  The components of $x_{\ell,k}$ are denoted $x_{\ell,k,i}$, for 
${i\in \ii{n_\ell}}$. We also denote
\[
G_{\ell,k} = \nabla_x^1f_{\ell}(x_{\ell,k}),
\ms
g_{\ell,k} \approx \nabla_x^1f_{\ell}(x_{\ell,k}),
\tim{  and }
H_{\ell,k} \approx \nabla_x^2f_{\ell}(x_{\ell,k}).
\]
If the bounds on the variables at level $\ell$ are $l_\ell$ and $u_\ell$, we have that $\calF_\ell = 
\prod_{i=1}^{n_\ell}[l_{\ell,i},u_{\ell,i}]$. In what follows, we assume that, for each $\ell$, 
$P_\ell$ and $R_\ell$ have non-negative entries. We define the vector of row-sums. 
$\sigma_\ell$, as given component-wise by
\beqn{sigma-def}
\sigma_{\ell,i} =\sum_{j \in \ii{n_{\ell-1}} | P_{\ell,i,j} >0 } P_{\ell,i,j}
= \sum_{j=1}^{n_{\ell-1}} P_{\ell,i,j} 
\ms ( i \in \ii{n_\ell} ).
\eeqn
When the $j$-th column is zero, we may choose $\sigma_{\ell,j}\geq0$ arbitrarily and 
we also define the maximum and minimum on the empty set 
$\{i\in\ii{n_\ell}|\mid P_{\ell,i,j} >0 \}$ to be are arbitrary with the latter 
not exceeding the former.
Using these conventions, we may now describe the \algname\ algorithm.

\algo{top-algo}{\algname}{
\begin{description}
\item[Initialization:] The bounds $l$ and $u$, the starting point $x_{r, 0}$, the hierarchy
  $\{f_\ell,\sigma_\ell, m_\ell\}_{\ell=0}^r$,
  $\{P_\ell,R_\ell\}_{\ell=0}^{r-1}$,
  and the constants $\kappa_s,\kap{2nd} \geq 1$
  and $\kap{1st},\kap{gs},\varsigma\in(0,1]$ are given.
\item[Top level optimization:]
return {\sf \algname-r}$(r,P_\calF(x_{r, 0}),l,u,\varsigma^2, 0, \infty, m_r)$
\end{description} 
}

\algo{rec-algo}{$x_{\ell,*} = $ {\sf \algname-r}
    ($\ell,x_{\ell,0},l_\ell,u_\ell,w_{\ell,-1}, \theta_{1,\ell},\theta_{2,\ell}, m_\ell$)}{
\begin{description}
\item[Step 0: Initialization:] Set $k=0$. 
\item[Step 1: Start iteration:]
  Compute $g_{\ell,k}$ as an approximation of $G_{\ell,k}$,
  \beqn{dk-def}
   d_{\ell,k} \eqdef P_{\calF_\ell}( x_{\ell,k} - g_{\ell,k} ) - x_{\ell,k},
  \eeqn
  \beqn{wDk-def}
  w_{\ell,k,i} = \sqrt{w_{\ell,k-1,i}^2+ d_{\ell,k,i}^2}
  \tim{and}
  \Delta_{\ell,k,i} =\frac{|d_{\ell,k,i}|}{w_{\ell,k,i}},
  \tim{ for } i\in \ii{n_\ell}.
  \eeqn
  If $\ell < r$ and $k = 0$, readjust
  \beqn{w0corr-def}
  w_{\ell,0} \leftarrow \max\left[1,
  \frac{\|\Delta_{\ell,0}\|}{\theta_{2,\ell}}\right]\,w_{\ell,0}
  \tim{ and }
  \Delta_{\ell,0} \leftarrow \min\left[1,
  \frac{\theta_{2,\ell}}{\|\Delta_{\ell,0}\|}\right]\,\Delta_{\ell,0}
  \eeqn
  \vspace*{-2mm}
  and return $x_{\ell,*}=x_{\ell,0}$ if
  \beqn{nogain}
  |d_{\ell,0}^T\Delta_{\ell,0}| < \theta_{1,\ell}.
  \eeqn
  Otherwise, set
  \beqn{Bk-def}
  \calB_{\ell,k}
  = \left\{ x \in \Re^{n_\ell} \mid | x_i - x_{\ell,k,i}| \leq \Delta_{\ell,k,i},
  \tim{for} i \in\ii{n_\ell} \right\}
  \eeqn
  \beqn{sL-def}
  s_{\ell,k}^L = P_{\calF_\ell\cap\calB_{\ell,k}}(x_{\ell,k}-g_{\ell,k})-x_{\ell,k} 
  \eeqn
  and select the type of iteration $k$ to be either `Taylor' or, if $\ell>0$, `recursive'.
\item[Step 2: Taylor iteration: ] 
  Choose $B_{\ell,k}$ a symmetric approximation of $H_{\ell,k}$ and compute
  \beqn{gamma-def}
  s_{\ell,k}^Q = \gamma_{\ell,k}s_{\ell,k}^L,
  \tim{where}
  \gamma_{\ell,k} = \left\{\begin{array}{ll}
      \min\left[ 1,\bigfrac{-g_{\ell,k}^Ts_{\ell,k}^L}{(s_{\ell,k}^L)^T B_{\ell,k} s_{\ell,k}^L}\right],
      & \tim{if } (s_{\ell,k}^L)^TB_{\ell,k}s_{\ell,k}^T > 0, \\
      1, & \tim{otherwise,}
    \end{array}\right.
  \eeqn
  select $s_{\ell,k}$ such that, for all $i\in\ii{n}$,
  \beqn{inside-and-decr}
  x_{\ell,k}+s_{\ell,k} \in \calF,
  \ms
  |s_{\ell,k,i}|\leq \kappa_s \Delta_{\ell,k,i}
  \tim{ and }
  g_{\ell,k}^Ts_{\ell,k}+\half s_{\ell,k}^TB s_{\ell,k}
  \leq \tau\left(g_{\ell,k}^Ts_{\ell,k}^Q+\half (s_{\ell,k}^Q)^TB
  s_{\ell,k}^Q\right),
  \eeqn
  \vspace*{-1mm}
  and go to Step~4.
\item[Step 3: Recursive iteration: ]
  Compute
     \beqn{kappadD-def}
     [ x_{\ell-1,0}, w_{\ell-1}]
     =  R_\ell [ x_{\ell,k}, w_{\ell,k}],
     \eeqn
     \beqn{lower-level-bounds}
\!\!\!\!
     l_{\ell-1,i} = x_{\ell-1,0,i} + \max_{q \mid P_{\ell,q,i}>0}\left[\frac{l_{\ell,q} - x_{\ell,k,q}}{\sigma_{\ell,q}}\right],
     \ms
     u_{\ell-1,i} = x_{\ell-1,0,i} + \min_{q \mid P_{\ell,q,i}>0} \left[\frac{u_{\ell,q} - x_{\ell,k,q}}{\sigma_{\ell,q}}\right],
     \eeqn
     \beqn{thetas-def}
     \theta_{1,\ell-1} = \kap{1st} \,|d_{\ell,k}^T\Delta_{\ell,k}|,
     \ms\ms
     \theta_{2,\ell-1} = \kap{2nd} \,\|s_{\ell,k}^L\|,
     \eeqn
     and set
  \vspace*{-1mm}
  \[
  \hspace*{-10mm}
  s_{\ell,k} = P_\ell \Big[ \mbox{\sf \algname-r}(\ell-1,x_{\ell-1,0},l_{\ell-1},u_{\ell-1},w_{\ell-1},\theta_{1,\ell-1}, \theta_{2,\ell-1}, m_{\ell-1}) -x_{\ell-1,0}\Big].
  \]
\item[Step 4: Loop:] If $\ell<r$ and
  \beqn{slope-cond}
   g_{\ell,0}^T(x_{\ell,k+1}-x_{\ell,0}) >\kappa_{gs}g_{\ell,0}^Ts_{\ell,0},
   \eeqn
   return $x_{\ell,*}=x_{\ell,k}$. Else set $x_{\ell,k+1} =
   x_{\ell,k}+s_{\ell,k}.$
   If $k = m_\ell$ return $x_{\ell,*}=x_{\ell,k+1}$.
   Otherwise, increment $k$ by one and go to Step~1.
\end{description} 
}
\newpage

\begin{enumerate}
\item The choice \( P_\ell = R_\ell^T \) corresponds to a Galerkin approach, commonly used in linear and nonlinear multilevel methods due to its symmetry and variational properties~\cite{trottenberg2001multigrid}. A more general choice leads to a Petrov-Galerkin formulation, which is particularly useful for problems that lack a variational structure, such as certain non-symmetric or indefinite problems, or problems arising in neural network training.

\item Observe that \req{nogain} is checked before any (gradient) evaluation at level $\ell$.
Should \req{nogain} hold, $x_{\ell,*}=x_{\ell,0}$ is returned, which corresponds to a
``void'' iteration (without any evaluation) at levels $\ell$ and $\ell+1$. In order to
simplify notations, we ignore such iterations below and assume that \req{nogain} fails
whenever {\sf \algname-r} is called at level $\ell<r$. Observe also
that
\beqn{mesure-1rst}
 |d_{\ell,k}^T\Delta_{\ell,k}| = |d_{\ell,k}|^T\Delta_{\ell,k} =\sum_{i=1}^{n_\ell} \frac{d_{\ell,k,i}^2}{w_{\ell,k,i}}
\eeqn
and hence \req{nogain} compares two non-negative numbers.
\item
The choice of $B_k$ in \req{gamma-def} is only restricted by the request that 
it should be bounded in norm (see AS.3 below). This allows for a wide range of 
approximations, such as Barzilai-Borwein or safeguarded (limited-memory) quasi-Newton.  
However, using the true Hessian $H_{\ell,k}$, that is computing 
$(s_{\ell,k}^L)^T H_{\ell,k} s_{\ell,k}^L$ instead of 
$(s_{\ell,k}^L)^T B_{\ell,k} s_{\ell,k}^L$, is numerically realistic 
(and often practically beneficial) even for very large problems if one uses 
(complex) finite-differences, thereby avoiding the need to evaluate 
full Hessians.
\item
As stated, the \algname\ algorithm has no stopping criterion.  In practice, termination 
of a particular realization may be decided if $\|d_{r,k}\| \leq
\epsilon$ in order to ensure an approximate $\epsilon$-first-order
critical point for problem \req{problem}.
\item
The point $x_{\ell,k} + s_{\ell,k}^Q$ minimizes the quadratic model of
the objective function's decrease given by
$g_{\ell,k}^Ts_{\ell,k}+\half s_{\ell,k}^TB s_{\ell,k}$ 
along $s_{\ell,k}$ and inside $\calB_{\ell,k}$.  In the terminology of
trust-region methods, it can therefore be considered as the ``Cauchy
point'' at iteration $(\ell,k)$. 
\item \label{diff-sL}
If only one level is considered ($r=1$), the algorithm is close to the
algorithm described in \cite{GratJeraToin24a} but yet differs from it
in a small but significant detail: the norm of the "linear" step
$\|s_{\ell,k}^L\|$ only results from \req{sL-def} and is independent
of the weights $w_{\ell,k}$. Thus, $|s_{\ell,k,i}|$ may be smaller
than $\Delta_{\ell,k,i}$, when $w_{\ell,k,i}<1$ and $\|s_{\ell,k}^L\|$
can  therefore be smaller than the distance from $x_{\ell,k}$ to the
boundary of $\calB_{\ell,k}$. This particular choice  makes the
restriction imposed by the second part of \req{inside-and-decr}
necessary for our multilevel argument. This latter restriction is
however only active when weights are smaller than one,  which is
typically only true for the first few iterations.  Thus \algname\  may
be slighly more cautious than the algorithm of \cite{GratJeraToin24a} close
to the starting point, which can be an advantage when early gradients 
may cause dangerously large steps.  However the inequality
\beqn{nsL}
\|s_{\ell,k}^L\| \leq \|\Delta_{\ell,k}\|
\eeqn
holds for all $\ell\in\iiz{r}$ and $k\geq 0$ because
$\calF \cap \calB_{\ell,k}\subseteq \calB_{\ell,k}$.
\item
The arguments $\theta_{1,\ell-1}$ and $\theta_{2,\ell-1}$ which are
passed to level $\ell-1$ in a recursive iteration are meant to
control, respectively, the minimal first-order achievement and the
maximal second-order effect of the step resulting from the recursive
call. This is obtained by imposing constraints on the zero-th
iteration at level $\ell-1$, in the sense that the mechanism of the
algorithm (see the second part of \req{wDk-def} and \req{thetas-def})
and our assumption that \req{nogain} fails) ensures the inequalities
\beqn{interlevels} 
|d_{\ell,0}^T\Delta_{\ell,0}| 
\geq \theta_{1,\ell} 
= \kap{1st} |d_{\ell+1,k}^T\Delta_{\ell+1,k}|
\tim{ and } 
\|s_{\ell,0}^L\|
\leq \|\Delta_{\ell,0}\| 
\leq \theta_{2,\ell} 
= \kap{2nd}\|s_{\ell+1,k}^L\|, 
\eeqn 
where iteration $(\ell+1,k)$ is the
upper-level ``parent'' iteration from which the call to {\sf \algname-r} at level
$\ell $ has been made. In addition, 
condition \req{slope-cond} is meant to guarantee 
that lower-level iterations $(\ell,k)$ for $k>0$ 
do not jeopardize the projected first-order progress achieved 
by iteration $(\ell,0)$. This condition is necessary because, again, 
nothing is known or assumed about the true or approximate gradients 
at lower-level iterations beyond the zero-th.
\end{enumerate}

To ensure that the \algname\ is well-defined, it remains to show that all iterates at each level remain feasible for the level-dependent bounds.  This is ensured by the following lemma.

\llem{lemma:feasibility}{
For each $\ell\in\iiz{r}$ and each $k\geq 0$, one has that $x_{\ell,k}\in \calF_\ell$.
}

\proof{
First note that $x_{r,0}\in\calF_r=\calF$  because
$x_{r,0}\in\calF$ by construction in the top level recursive call in \algname.
Suppose now that $x_{\ell,k}\in \calF_\ell$ for some $\ell\in\iiz{r}$ and some $k\geq0$
and consider the next computed iterate.  Three cases are possible.
\begin{itemize}
\item The first is when iteration $(\ell,k)$ is a Taylor iteration.  In this case, the first part of \req{inside-and-decr} ensures that $x_{\ell,k+1}= x_{\ell,k}+s_{\ell,k}\in \calF_\ell$, as requested.
\item The second is when iteration $(\ell,k)$ is recursive and the next iterate is $x_{\ell-1,0}$.
If this is the case, \req{lower-level-bounds} ensures that $x_{\ell-1,0}\in\calF_{\ell-1}$.
\item The third is when $x_{\ell,k}$ is the final iterate $x_{\ell,*}$ at level $\ell<r$ 
of a  sequence of step initiated by a recursive iteration $(\ell+1,p)$, in which case the 
next iterate is 
\beqn{case3}
x_{\ell+1,p+1} = x_{\ell+1,p} + s_{\ell+1,p}
= x_{\ell+1,p}+ P_{\ell+1} \widehat{s}_{\ell,k}.
\tim{ where } \widehat{s}_{\ell,k} = x_{\ell,k}-x_{\ell,0}.
\eeqn
Because $x_{\ell,k}\in \calF_\ell$, we have that, for all $j\in\ii{n_\ell}$,
\beqn{low-level-feasible}
l_{\ell,j}
\leq x_{\ell,0,j} + \widehat{s}_{\ell,k,j}
\leq u_{\ell,j}.
\eeqn
Now define, for each $i\in\ii{n_{\ell+1}}$ and $j\in\ii{n_\ell}$,
\[
\calR_{\ell+1,i} = \{ j \in \ii{n_\ell} \mid P_{\ell+1,i,j} > 0 \}
\tim{ and }
\calC_{\ell+1,j} = \{ i \in \ii{n_\ell+1}\mid P_{\ell+1,i,j} > 0 \},
\]
the supports of row $i$ and column $j$ of $P_{\ell+1}$, respectively.
Thus 
\beqn{Pshat}
[P_{\ell+1} \widehat{s}_{\ell,k}]_i= \bigsum_{j\in\calR_{\ell+1,i}}P_{\ell+1,i,j}\,\widehat{s}_{\ell,k,j},
\eeqn
for $i\in \ii{n_{\ell+1}}$ because we assumed that
all entries of $P_{\ell+1}$ are non-negative.
Observe also that, if $j\in\calR_{\ell+1,i}$, then $P_{\ell+1,i,j} > 0$,
$i\in\calC_{\ell+1,j}\neq \emptyset$ and $\sigma_{\ell+1,i}>0$.  Thus, because of \req{lower-level-bounds} and \req{low-level-feasible},
\[
\frac{l_{\ell+1,i}-x_{\ell+1,p,i}}{\sigma_{\ell+1,i}}
\leq \max_{q\in\calC_{\ell+1,j}}\left[\frac{l_{\ell+1,q}-x_{\ell+1,p,q}}{\sigma_{\ell+1,q}}\right]
\leq \widehat{s}_{\ell,k,j}, 
\]
and therefore, remembering that $l_{\ell+1,i}-x_{\ell+1,p,i}\leq 0$ and using \req{sigma-def},
\beqn{lower-ok}
\sum_{j\in\calR_{\ell+1,i}}P_{\ell+1,i,j}\,\widehat{s}_{\ell,k,j}
\geq \frac{l_{\ell+1,i}-x_{\ell+1,p,i}}{\sigma_{\ell+1,i}}\bigsum_{j\in\calR_{\ell+1,i}}P_{\ell+1,i,j}
= l_{\ell+1,i}-x_{\ell+1,p,i},
\eeqn
for $i\in\ii{n_{\ell+1}}$. 
Similarly, for $j\in\calR_{\ell+1,i}$,
\[
\widehat{s}_{\ell,k,j} 
\leq \min_{q\in\calC_{\ell+1,j}}\left[\frac{u_{\ell+1,q}-x_{\ell+1,p,q}}{\sigma_{\ell+1,q}}\right]
\leq \frac{u_{\ell+1,i}-x_{\ell+1,p,i}}{\sigma_{\ell+1,i}},
\] 
and thus, for $i\in\ii{n_{\ell+1}}$,
\beqn{upper-ok}
\bigsum_{j\in\calR_{\ell+1,i}}P_{\ell+1,i,j}\,\widehat{s}_{\ell,k,j} 
\leq \bigfrac{u_{\ell+1,i}-x_{\ell+1,p,i}}{\sigma_{\ell+1,i}}\bigsum_{j\in\calR_{\ell+1,i}}P_{\ell+1,i,j}
= u_{\ell+1,i}-x_{\ell+1,p,i}.
\eeqn
Combining \req{Pshat}, \req{upper-ok} and \req{lower-ok} gives that
\[
l_{\ell+1,i}-x_{\ell+1,p,i}
\leq [P_{\ell+1} \widehat{s}_{\ell,k}]_i
\leq u_{\ell+1,i}-x_{\ell+1,p,i},
\]
for each $i\in \ii{n_{\ell+1}}$ and therefore, using \req{case3}, that $x_{\ell+1,p+1}\in\calF_{\ell+1}$.
\end{itemize}
As a consequence, we see that feasibility with respect to the relevant level-dependent feasible set is maintained at all steps of the computation, yielding the desired result.
}

\section{Convergence analysis}
\label{section:convergence}

The convergence analysis of the
\algname\ algorithm can be viewed as an extension of that proposed, in
the single-level case, for the \al{ADAGB2} algorithm in 
\cite{BellGratMoriToin25}. We restate here the results of this reference 
that are necessary for our new argument, but the reader is referred to 
\cite{BellGratMoriToin25} for their proofs, which are based on a 
subset of the following assumptions.\\ 

\noindent
{\bf AS.1:} \textit{Each $f_\ell$ for $\ell\in\iiz{r}$ is twice continuously differentiable.}\\

\noindent
{\bf AS.2:} \textit{For each $\ell\in\iiz{r}$ there exists a constant $L_\ell \geq0$ such that for all $x,y \in \Re^{n_\ell}$
\[
\|\nabla_x^1f_\ell(x)-\nabla_x^1f_\ell(y)\| \leq L_\ell\|x-y\|.
\]
}\\

\noindent
{\bf AS.3:} \textit{There exists a constant $\kappa_B\geq 1$ such that, for each $\ell\in\iiz{r}$
and all $k \geq 0$, $\|B_{\ell,k}\| \leq \kappa_B$.}\\

\noindent
{\bf AS.4:} \textit{The objective function is bounded below on the feasible
domain, that is there exists a constant $\flow < f_r(x_{r,0})$, such that $f_r(x)\geq
\flow$ for every $x \in \calF$.}\\

In order to state our remaining assumptions, we define the event
\beqn{A-def}
\calE = \Big\{ \|d_{r,0}\|^2 \geq \varsigma \Big\}.
\eeqn
This event occurs or does not occur at iteration $(r,0)$, i.e.\ at the very
beginning of a realization of the \algname\ algorithm.  The convergence theory which
follows is dependent on the (in practice extremely likely) occurrence of $\calE$ and our subsequent
stochastic assumptions are therefore conditioned by this event.  They will be
formally specified by considering, at iteration
$(\ell,j)$, expectations conditioned by the past iterations and by $\calE$, which will
be denoted by the symbol $\EcondA{\ell,j}{\cdot}$. Note that,
because $\calE$ is measurable for all iterations after $(r,0)$,
$\EcondA{\ell,j}{\cdot}=\mathbb{E}_{\ell,j}\![\cdot]$, whenever $\ell < r$ or $k>0$.\\

\noindent
{\bf AS.5:}
\textit{There exists a constant $\kappa_{Gg}\geq 0$, such that
\beqn{AS6top}
\EcondA{\ell,k}{\|g_{\ell,k}-G_{\ell,k}\|^2} \leq \kappa_{Gg}^2
\,\EcondA{\ell,k}{\min[\|s_{\ell,k}\|^2,\theta_{2,\ell}^2]},
\eeqn
for all iterations $(\ell,k)$.}\\

\noindent
{\bf AS.6:}
\textit{There exists a constant $\kappa_\tau>0$, such that
\beqn{coherence-recursive}
\EcondA{\ell,0}{\|P_{\ell+1}^T G_{\ell+1,k} - g_{\ell,0}\|^2}
\leq \kappa_\tau^2 \EcondA{\ell,0}{\theta_{2,\ell}^2},
\eeqn
for all recursive iterations $(\ell+1,k)$.}\\

\noindent
Note that, since
$
\EcondA{\ell,k}{\min[\|s_{\ell,k}\|^2,\theta_{2,\ell}^2]}
\leq \min\left[\EcondA{\ell,k}{\|s_{\ell,k}\|^2},
  \EcondA{\ell,k}{\theta_{2,\ell}^2}\right],
$
AS.5 implies that
\beqn{AS6b}
\EcondA{\ell,k}{\|g_{\ell,k}-G_{\ell,k}\|^2} \leq \kappa_{Gg}^2
 \min\left[\EcondA{\ell,k}{\|s_{\ell,k}\|^2},
  \EcondA{\ell,k}{\theta_{2,\ell}^2}\right],
\eeqn
for all $(\ell,k)$. Note that AS.5 implies the directional error condition of \cite[AS.5]{BellGratMoriToin25},
while AS.6 is specific to the multi-level context. It is the only 
assumption relating $f_\ell(\cdot)$ and its model
$f_{\ell-1}(\cdot)$, and is limited to requiring a weak probabilistic
coherence between gradients of these two functions at $x_{\ell+1,k}$
and its restriction $x_{\ell,0} = R_\ell x_{\ell+1,k}$. Nothing is
assumed for the approximate gradients at iterates $x_{\ell,k}$ for
$k>0$.  This makes the choice of lower-level models very open.

The following ``linear descent'' lemma is a stochastic and 
bound-constrained generalization of \cite[Lemma~2.1]{GratJeraToin24a} 
(which, in the unconstrained and deterministic context, uses a 
different optimality measure and a different
definition of $s_{\ell,k}^L$).

\llem{lemma:single-level-decrease}{
Suppose that AS.3 and AS.5 hold and that iteration $(\ell,k)$ is a
Taylor iteration. Then
\beqn{gen-decr}
\EcondA{\ell,k}{G_{\ell,k}^Ts_{\ell,k}}
\leq
-\frac{\tau\varsigma^2}{2\kappa_B}\EcondA{\ell,k}{\,|d_{\ell,k}^T\Delta_{\ell,k}|\,}
+ \kappa_s^2(\half\kappa_B+\kappa_{Gg})\EcondA{\ell,k}{\|\Delta_{\ell,k}\|^2}.
\eeqn
}

\proof{See \cite[Lemma~2.1]{BellGratMoriToin25}.}

\noindent
If there is only one level and given \req{gen-decr}, the argument for 
proving convergence (albeit not complexity) is now relatively intuitive. 
One clearly sees in this relation that the first-order Taylor 
approximation decreases because of the first term on the right-hand 
but possibly increases because of the second-order effects of the second term.
When the weights $w_{\ell,j,i}$ become large, which, because of the
first part of \req{wDk-def}, must happen if convergence
stalls, then the first-order terms eventually dominate, causing a
decrease in the objective function, in turn leading to convergence. 
Because we have used the boundedness of $B_k$ only in the direction $s_{\ell,k}^L$,
we note that AS.3 could be weakened to merely require that 
${(s_{\ell,k}^L)^TB_{\ell,k}s_{\ell,k}^L}\leq \kappa_B \|s_{\ell,k}^L\|^2$ for 
all $s_{\ell,k}^L$ generated by the algorithm.

The challenge for proving convergence in  the multilevel context, 
which is the focus of this paper, is to handle, at level 
$\ell+1$, the effect of iterations at lower levels.  While AS.6 
suggests that the first iteration at level $\ell$ is somehow 
linked to the information available at iteration $k$ of level 
$\ell+1$, further iterations at level $\ell$ are only constrained 
by the condition \req{slope-cond}.  This condition enforces a 
minimum descent up to first order, but does not limit 
$\|x_{\ell,*}-x_{\ell,0}\|$, the length of the lower level step  
obtained at a recursive iteration (not imposing such a 
limitations appeared to be computationally advantageous).  
However, the convergence argument using \req{gen-decr} we have just outlined indicates that second-order effects 
(depending on this length) must be controlled compared to first-
order ones.  Thus it is not surprising that our first technical 
lemma considers the length of recursive steps. It is proved using 
the useful inequality 
\beqn{norm-sum-squared}
\|\sum_{i=1}^q a_i \|^2 \leq q \sum_{i=1}^q\| a_i \|^2,
\eeqn
which is valid for any collection $\{a_i\}_{i=1}^q$ of vectors.

\llem{lemma:length-recursive-step}{
Suppose that AS.1, AS.2, AS.5 and AS.6 hold. Suppose also that 
iteration $(\ell+1,k)$ is a recursive iteration whose lower-level 
iterate is $x_{\ell,*}=x_{\ell,p}$ for some 
$p \in \iiz{m_\ell+1}$. Then, for all $j\in\ii{p}$,
\beqn{xjmx0-bound}
\EcondA{\ell+1,k}{\|x_{\ell,j}-x_{\ell,0}\|^2}
\leq \alpha_\ell\,\EcondA{\ell+1,k}{\theta_{2,\ell}^2}
\eeqn
where
\beqn{alpha1-def}
\alpha_\ell \eqdef \left(\bigfrac{96}{\varsigma} \kappa_{Gg}^2+4\right)
\sum_{t=0}^{m_\ell}\left[\bigfrac{4\kappa_s^2}{\varsigma}(12 L_\ell^2+6)+2\right]^t.
\eeqn
Moreover, for all $j\in\iiz{p-1}$
\beqn{nsj-bound}
\EcondA{\ell+1,k}{\|s_{\ell,j}\|^2}
\leq 4 \alpha_\ell\, \EcondA{\ell+1,k}{\theta_{2,\ell}^2}.
\eeqn
}

\proof{
Consider  $j\in\ii{p}$. We start by
observing that, using \req{norm-sum-squared}, AS.2 and \req{AS6b}, 
\beqn{gamma}
\begin{array}{lcl}
\EcondA{\ell,j}{\|g_{\ell,j}-g_{\ell,0}\|^2}
& \leq & 3\Big(\EcondA{\ell,j}{\|g_{\ell,j}-G_{\ell,j}\|^2}
               +\EcondA{\ell,j}{\|G_{\ell,j}-G_{\ell,0}\|^2}\\
&& ~~          +\EcondA{\ell,0}{\|g_{\ell,0}-G_{\ell,0}\|^2}\Big)\\
& \leq &  3\left(L_\ell^2\EcondA{\ell,j}{\|x_{\ell,j}-x_{\ell,0}\|^2}
               + 2\kappa_{Gg}^2\EcondA{\ell,0}{\theta_{2,\ell}^2} \right).
\end{array}
\eeqn
Using \req{norm-sum-squared}, \req{sL-def} and the contractivity of the projection $P_\calF$, we also obtain that
\beqn{delta}
\begin{array}{lcl}
\EcondA{\ell,j}{\|d_{\ell,j}-d_{\ell,0}\|^2}
& = & \EcondA{\ell,j}{\|P_\calF(x_{\ell,j}-g_{\ell,j})-x_{\ell,j}-P_\calF(x_{\ell,0}-g_{\ell,0})+x_{\ell,0}\|^2}\\*[1ex]
& \leq & 2\, \EcondA{\ell,j}{\|P_\calF(x_{\ell,j}-g_{\ell,j})-P_\calF(x_{\ell,0}-g_{\ell,0})\|^2}
         + 2\,\EcondA{\ell,j}{\|x_{\ell,j}-x_{\ell,0}\|^2}\\*[1ex]
& \leq &2\,\EcondA{\ell,j}{\| (x_{\ell,j}-g_{\ell,j})-(x_{\ell,0}-g_{\ell,0})\|^2}
         + 2\,\EcondA{\ell,j}{\|x_{\ell,j}-x_{\ell,0}\|^2}\\*[1ex]
& \leq &4\,\EcondA{\ell,j}{\|g_{\ell,j}-g_{\ell,0}\|^2}
         + 6\, \EcondA{\ell,j}{\|x_{\ell,j}-x_{\ell,0}\|^2}.\\*[1ex]
\end{array}
\eeqn
Substituting \req{gamma} into \req{delta} gives that
\beqn{dd2}
\EcondA{\ell,j}{\|d_{\ell,j}-d_{\ell,0}\|^2}
\leq (12 L_\ell^2+6)\EcondA{\ell,j}{\|x_{\ell,j}-x_{\ell,0}\|^2}
   + 24 \kappa_{Gg}^2 \EcondA{\ell,0}{\theta_{2,\ell}^2}.
\eeqn
But, for each $i\in\ii{n_{\ell}}$, the nondecreasing nature of
$w_{\ell,j-1,i}$ as a function of $j$ and the triangle inequality give that
\[
s_{\ell,j-1,i}^L
\leq \frac{|d_{\ell,j-1,i}|}{w_{\ell,j-1,i}}
\leq \left|\frac{d_{\ell,j-1,i}}{w_{\ell,0,i}}\right|
\leq \left|\frac{d_{\ell,j-1,i}-d_{\ell,0,i}}{w_{\ell,0,i}} + \frac{d_{\ell,0,i}}{w_{\ell,0,i}}\right|.
\]
Using \req{norm-sum-squared} once more with \req{wDk-def}, we obtain that
\[
(s_{\ell,j-1,i}^L)^2
\leq  2\left|\frac{d_{\ell,j-1,i}-d_{\ell,0, i}}{w_{\ell,0,i}}\right|^2
     +2\left|\frac{d_{\ell,0,i}}{w_{\ell,0,i}}\right|^2
\leq  \frac{2}{\varsigma}|d_{\ell,j-1,i}-d_{\ell,0, i}|^2
     +2\,\Delta_{\ell,0,i}^2
     \]
and thus
\[
\begin{array}{lcl}
\EcondA{\ell,j-1}{\|s_{\ell,j-1,i}^L\|^2}
& \leq & \bigfrac{2}{\varsigma}\,\EcondA{\ell,j-1}{\|d_{\ell,j-1,i}-d_{\ell,0, i}\|^2}
    +2\,\EcondA{\ell,j-1}{\|\Delta_{\ell,0}\|^2}.
\end{array}
\]
Substituting now \req{dd2} for $j-1$ in this inequality, using the second part of
\req{interlevels} and the tower property, we deduce that, for $j\in\iibe{2}{p}$,
\beqn{bigs}
\EcondA{\ell+1,k}{\|s_{\ell,j-1,i}^L\|^2}
\leq \bigfrac{2}{\varsigma}(12 L_\ell^2+6)\,\EcondA{\ell+1,k}{\|x_{\ell,j-1}-x_{\ell,0}\|^2}
+ \left(\bigfrac{48}{\varsigma}\kappa_{Gg}^2+2\right) \,\EcondA{\ell+1,k}{\theta_{2,\ell}^2}.
\eeqn
Hence, using \req{norm-sum-squared}, this last inequality and the second part of \req{interlevels}, we
obtain that, for ${j\in\iibe{2}{p}}$,
\beqn{xx}
\begin{array}{lcl}
\EcondA{\ell+1,k}{\|x_{\ell,j}-x_{\ell,0}\|^2}
& \leq & 2\,\EcondA{\ell+1,k}{\|s_{\ell,j-1}\|^2} + 2\,\EcondA{\ell+1,k}{\|x_{\ell,j-1}-x_{\ell,0}\|^2},\\*[2ex]
& \leq & \left[\bigfrac{4\kappa_s^2}{\varsigma}(12 L_\ell^2+6)+2\right]\,\EcondA{\ell+1,k}{\|x_{\ell,j-1}-x_{\ell,0}\|^2},\\*[2ex]
&& ~~ + \left(\bigfrac{96}{\varsigma}\kappa_{Gg}^2+4\right) \,\EcondA{\ell+1,k}{\theta_{2,\ell}^2}.
\end{array}
\eeqn
Applying this inequality recursively down to $j=0$ then gives that
\[
\EcondA{\ell+1,k}{\|x_{\ell,j}-x_{\ell,0}\|^2}
\leq \left[\left(\bigfrac{96}{\varsigma}\kappa_{Gg}^2+4\right)
\sum_{t=0}^{j-1}\left[\bigfrac{4\kappa_s^2}{\varsigma}(12 L_\ell^2+6)+2\right]^t\right]
\,\EcondA{\ell+1,k}{\theta_{2,\ell}^2},
\]
which, with the bound $j \leq p \leq m_\ell+1$, gives \req{xjmx0-bound} and \req{alpha1-def}.
Moreover \req{xjmx0-bound} and
\[
\EcondA{\ell,j}{\|s_{\ell,j}\|^2}
\leq  2\left(\EcondA{\ell,j}{\|x_{\ell,j+1}-x_{\ell,0}\|^2}
             +\EcondA{\ell,j}{\|x_{\ell,j}-x_{\ell,0}\|^2}\right)
\]
imply \req{nsj-bound} for $j\leq p-1$.
}

\noindent
Once an upper bound on the length of the total lower-level step (given
by \req{xjmx0-bound}) is known, we may use it to consider the relation
between the expected linear decrease at level $\ell+1$ given knowledge
of its value at level $\ell$.

\llem{lemma:recursive-decrease}{
Suppose that AS.1, AS.2, AS.5 and AS.6 hold. Suppose also that
iteration $(\ell+1,k)$ is a recursive iteration whose final low-level
iterate is $x_{\ell,*}=x_{\ell,p}$ for some $p \in \iiz{m_\ell+1}$.
Suppose furthermore that, for some constants
$\beta_{1,\ell},\beta_{2,\ell} >0$, the first iteration $(\ell,0)$ at
level $\ell$ is such that
\beqn{lowlvl-cond}
\EcondA{\ell+1,k}{G_{\ell,0}^Ts_{\ell,0}}
\leq - \beta_{1,\ell} \EcondA{\ell+1,k}{\,|d_{\ell,0}^T\Delta_{\ell,0}|\,}
  + \beta_{2,\ell} \EcondA{\ell+1,k}{\|s_{\ell,0}^L\|^2}.
\eeqn
Then
\beqn{recursive-linear-decrease}
\EcondA{\ell+1,k}{G_{\ell+1,k}^Ts_{\ell+1,k}}
\leq -\beta_{1,\ell+1}\EcondA{\ell+1,k}{\,|d_{\ell+1,k+1}^T\Delta_{\ell+1,k+1}|\,}
+\beta_{2,\ell+1}\EcondA{\ell+1,k}{\|s_{\ell+1,k}^L\|^2},
\eeqn
where 
\beqn{beta-rec12}
\beta_{1,\ell+1} = \kap{1st} \beta_{1,\ell},
\ms
\beta_{2,\ell+1} =
\kap{2nd}^2(\kappa_{gs}\beta_{2,\ell}+\kappa_\tau\bigsqrt{\alpha_\ell}+\kappa_s\kappa_{Gg})
\eeqn
with $\alpha_\ell$ defined by \req{alpha1-def}.
}

\proof{
We first observe that \req{slope-cond} in Step~4 of the algorithm
ensures that
\[
\EcondA{\ell+1,k}{g_{\ell,0}^T(x_{\ell,p}-x_{\ell,0})}
\leq \kappa_{gs}\,\EcondA{\ell+1,k}{g_{\ell,0}^Ts_{\ell,0}}.
\]
Thus, we have that
\beqn{eq3}
\begin{array}{lcl}
\EcondA{\ell+1,k}{G_{\ell+1,k}^Ts_{\ell+1,k}}  
& = & \EcondA{\ell+1,k}{G_{\ell+1,k}^TP_{\ell+1}(x_{\ell,p}-x_{\ell,0})}\\
& = & \EcondA{\ell+1,k}{\Big(P_{\ell+1}^T G_{\ell+1,k}\Big)^T(x_{\ell,p}-x_{\ell,0}) }\\*[1ex]
& \leq & \EcondA{\ell+1,k}{g_{\ell,0}^T(x_{\ell,p}-x_{\ell,0}) }
    + |\EcondA{\ell+1,k}{(P_{\ell+1}^T G_{\ell+1,k} - g_{\ell,0})^T(x_{\ell,p}-x_{\ell,0})}|\\*[2ex]
& \leq & \kappa_{gs}\,\EcondA{\ell+1,k}{g_{\ell,0}^Ts_{\ell,0}}\\*[2ex]
&&~~    + \bigsqrt{\EcondA{\ell+1,k}{\|P_{\ell+1}^T G_{\ell+1,k} -
        g_{\ell,0}\|^2}\,\EcondA{\ell+1,k}{\|x_{\ell,p}-x_{\ell,0}\|^2}}.\\*[2ex]
\end{array}
\eeqn
In order to bound the term
$\EcondA{\ell+1,k}{\|x_{\ell,p}-x_{\ell,0}\|^2}$ in the right-hand
side of this inequality, we may now apply
Lemma~\ref{lemma:length-recursive-step} to iteration $(\ell,j)$.
Using the tower property, AS.6 and \req{xjmx0-bound}, we thus obtain that
\beqn{eq2}
\begin{array}{l}
\bigsqrt{\EcondA{\ell+1,k}{\|P_{\ell+1}^T G_{\ell+1,k} -
        g_{\ell,0}\|^2}\,\EcondA{\ell+1,k}{\|x_{\ell,p}-x_{\ell,0}\|^2}}\\*[2ex]
\hspace*{20mm} = \bigsqrt{\EcondA{\ell+1,k}{\EcondA{\ell,0}{\|P_{\ell+1}^T G_{\ell+1,k} -
          g_{\ell,0}\|^2}}\,\EcondA{\ell+1,k}{\|x_{\ell,p}-x_{\ell,0}\|^2}},\\*[2ex]
\hspace*{20mm} = \kappa_\tau\bigsqrt{\EcondA{\ell+1,k}{\EcondA{\ell,0}{\theta_{2,\ell}^2}}\,
         \alpha_\ell\,\EcondA{\ell+1,k}{\theta_{2,\ell}^2}},\\*[2ex]
\hspace*{20mm} = \kappa_\tau\bigsqrt{\alpha_\ell}\,\EcondA{\ell+1,k}{\theta_{2,\ell}^2}.
\end{array}
\eeqn
Substituting \req{eq2} in \req{eq3} and using \req{AS6b} and the
tower property then gives that
\[
\begin{array}{lcl}
\EcondA{\ell+1,k}{G_{\ell+1,k}^Ts_{\ell+1,k}}    
& \leq & \kappa_{gs}\, \EcondA{\ell+1,k}{g_{\ell,0}^Ts_{\ell,0}}
    + \kappa_\tau\bigsqrt{\alpha_\ell}\,\EcondA{\ell+1,k}{\theta_{2,\ell}^2}\\*[2ex]
& = & \kappa_{gs}\,\EcondA{\ell+1,k}{\EcondA{\ell,0}{G_{\ell,0}^Ts_{\ell,0}}}
    +\EcondA{\ell+1,k}{\EcondA{\ell,0}{(g_{\ell,0}-G_{\ell,0})^Ts_{\ell,0}}}\\*[2ex]
&&~~+\kappa_\tau\bigsqrt{\alpha_\ell}
    \,\EcondA{\ell+1,k}{\theta_{2,\ell}^2}\\*[2ex]
& = & \kappa_{gs}\,\EcondA{\ell+1,k}{\EcondA{\ell,0}{G_{\ell,0}^Ts_{\ell,0}}}
    +\EcondA{\ell+1,k}{\sqrt{\EcondA{\ell,0}{\|g_{\ell,0}-G_{\ell,0}\|^2}\EcondA{\ell,0}{\|s_{\ell,0}\|^2}}}\\*[2ex]
&&~~+\kappa_\tau\bigsqrt{\alpha_\ell}\,\EcondA{\ell+1,k}{\theta_{2,\ell}^2}\\*[2ex]
& = & \kappa_{gs}\,\EcondA{\ell+1,k}{\EcondA{\ell,0}{G_{\ell,0}^Ts_{\ell,0}}}
    +\kappa_s\kappa_{Gg}\EcondA{\ell+1,k}{\EcondA{\ell,0}{\theta_{2,\ell}^2}}\\*[2ex]
&&~~+\kappa_\tau\bigsqrt{\alpha_\ell}
    \,\EcondA{\ell+1,k}{\theta_{2,\ell}^2}\\*[2ex]
& = & \kappa_{gs}\,\EcondA{\ell+1,k}{G_{\ell,0}^Ts_{\ell,0}}
  +(\kappa_\tau\bigsqrt{\alpha_\ell}^2+\kappa_s\kappa_{Gg})
     \,\EcondA{\ell+1,k}{\theta_{2,\ell}^2}.\\*[2ex]
\end{array}
\]
We may now recall our recurrence assumption on the zero-th iteration
at level $\ell$ by applying \req{lowlvl-cond} and derive, using
\req{interlevels}, that  
\[
\begin{array}{lcl}
\EcondA{\ell+1,k}{G_{\ell+1,k}^Ts_{\ell+1,k}}
&\leq & - \kappa_{gs}\, \beta_{1,\ell} \EcondA{\ell+1,k}{\,|d_{\ell,0}^T\Delta_{\ell,0}|\,}
     + \kappa_{gs} \beta_{2,\ell} \EcondA{\ell+1,k}{\|s_{\ell,0}^L\|^2}\\*[1ex]
&&~~ + (\kappa_\tau\bigsqrt{\alpha_\ell}+\kappa_{Gg}^2)
\,\EcondA{\ell+1,k}{\theta_{2,\ell}^2}\\*[1ex]
&\leq & - \kappa_{gs}\beta_{1,\ell} \kap{1st}\EcondA{\ell+1,k}{\,|d_{\ell+1,k}^T\Delta_{\ell+1,k}|\,}\\*[1ex]
&&~~ + \kap{2nd}^2(\kappa_{gs}\beta_{2,\ell}+\kappa_\tau\bigsqrt{\alpha_\ell}+\kappa_s\kappa_{Gg})
  \,\EcondA{\ell+1,k}{\|s_{\ell+1,k}^L\|^2}.
\end{array}
\]
This proves \req{recursive-linear-decrease} with \req{beta-rec12}.
} 

\noindent
This lemma tells up how first- and second-order quantities behave when
one moves one level up. We now use these bounds to analyze the
complete hierarchy for $\ell=0$ up to $\ell = r$ and derive a
bound on the expectation of the linear decrease in the multi-level case. 

\llem{lemma:the-decrease}{
Suppose that AS.1, AS.2, AS.5 and AS.6 hold and consider an iteration $(r,k)$ at the top level.  Then
\beqn{the-decrease}
\EcondA{r,k}{G_{r,k}^Ts_{r,k}}
\leq -\beta_{r,1}\EcondA{r,k}{\,|d_{r,k}^T\Delta_{r,k}|\,}
+ \beta_{2,r} \EcondA{r,k}{\|\Delta_{r,k}\|^2},
\eeqn
for some constants $\beta_{r,1}$ and $\beta_{2,r}$ only dependent on 
the problem and algorithmic constants.
}

\proof{
First note that, for any Taylor iteration $(\ell,k)$, the second part of 
\req{inside-and-decr} and \req{gen-decr} ensures that \req{lowlvl-cond} holds with 
\beqn{betas-Taylor}
\beta_{1,\ell} = \frac{\tau\varsigma^2}{2\kappa_B}
\tim{ and }
\beta_{2,\ell} = \kappa_s^2(\half \kappa_B+\kappa_{Gg}).
\eeqn
Suppose now that level $\ell$ is the deepest level reached during the 
computation of the step $s_{r,k}$ at the top level.  By construction, 
all iterations at level $\ell$ are Taylor iterations and thus satisfy 
the condition \req{lowlvl-cond}. 
As a consequence, we may apply Lemma~\ref{lemma:recursive-decrease} 
and deduce that \req{recursive-linear-decrease} holds for the ``parent'' iteration $(\ell+1,k)$ from which the
recursion to level $\ell$ occurred. 
Thus all iterations (Taylor or recursion) at level $\ell+1$ in turn satisfy 
condition \req{lowlvl-cond} (with $\ell$ replaced by $\ell+1$) and
updated coefficients given by  \req{beta-rec12}.
Now the inequality $\kap{1st}\leq 1$  and  the relations \req{beta-rec12}  
ensure that, for all $t\in\iibe{\ell}{r-1}$,
\beqn{betas-rec}
\beta_{1,t+1} \leq \beta_{1,t},
\ms
\beta_{2,t+1} \geq \beta_{2,t}.
\eeqn
As a consequence, we may pursue the recurrence \req{beta-rec12} 
until $\ell+1 = r$, absorbing, if any, Taylor iterations and less deep
recursive ones, to finally define $\beta_{1,r}$ and $\beta_{2,r}$.  The
desired conclusion then follows from \req{recursive-linear-decrease}
and the second part of \req{interlevels}.
} 

\noindent
Comparing \req{gen-decr} and \req{the-decrease}, we observe that the
first-order behaviour of the multilevel algorithm is identical to that of the single-level version of the algorithm
described in \cite{BellGratMoriToin25}, except that the single level constants
\req{betas-Taylor} are now replaced by $\beta_{1,r}$ and $\beta_{2,r}$
as resulting from Lemma~\ref{lemma:the-decrease}.  It is
therefore not surprising that the line of argument used in
\cite[Theorem~3.2]{GratJeraToin22a} for the deterministic  unconstrained case and extended to the
stochastic bound-constrained case in \cite{BellGratMoriToin25}  can be invoked  to
cover the more general stochastic bound-constrained multi-level context.

\lthm{theorem:convergence}{
Suppose that AS.1--AS.6 hold and that the
\algname\ algorithm is applied to problem \req{problem}. Then
\beqn{gradbound}
\EA{\average_{j\in\iiz{k}}\|d_{r,j}\|^2}
\le \frac{\kap{conv}}{k+1},
\eeqn
with 
\beqn{kapconv-def}
\kap{conv}
=\bigfrac{\varsigma\,\kappa_W^2}{2}\,\left|W_{-1}\left(-\frac{1}{\kappa_W}\right)\right|^2
\leq \bigfrac{\varsigma\,\kappa_W^2}{2}\left|\log(\kappa_W) + \sqrt{2(\log(\kappa_w)-1)}\right|^2,
\eeqn
where $W_{-1}$ is the second branch of the Lambert function, $\Gamma_0 \eqdef f_r(x_0)-\flow$ and 
\beqn{kappaW-def}
\kappa_W =  \frac{4}{\beta_{1,r}\sqrt{\varsigma}}\max[1,\Gamma_0]
      \max\left[2,\bigfrac{n_r}{\Gamma_0}\max[ 1,\beta_{2,r} + \half \kappa_s^2 L_r]\right]
\eeqn
with $\beta_{r,1}$ and $\beta_{2,r}$ constructed using the recurrence \req{beta-rec12}.
}

\noindent
We may now return to the comment made after Lemma~\ref{lemma:single-level-decrease}
on the balance of first-order and second-order terms.  In the recurrence 
\req{beta-rec12}, from which the constants in Theorem~\ref{theorem:convergence} 
are deduced, the first -order terms (enforcing descent for large weights) are 
multiplied by $\beta_{1,r}$ which, for reasonable values of $\kappa_\delta$, 
does not decrease too quickly.  By contrast, the second-order terms involve 
the potentially large $\beta_{2,r}$.  From the theoretical point of view, 
it would therefore be advantageous to limit the size of recursive steps to a 
multiple of $|d_{r,k,i}|/w_{r,k,i}$, as was done in \cite{GratKopaToin23}. 
However, this restriction turned out to be counter-productive in practice, 
which prompted the considerably more permissive approach developed here.  
Unsurprisingly, letting lower-level iterations move far from $x_{\ell,0}$ 
comes at the cost of repeatedly using the Lipschitz assumption to estimate 
the distance of the lower-level iterates to $x_{\ell,0}$, (see 
\req{xjmx0-bound}) which causes the potentially large value of $\beta_{2,r}$ 
and, consequently, of $\kappa_W$. 

There are however ways to improve the constants of
Theorem~\ref{theorem:convergence}.
As we suggested above, the most obvious way to improve on second-order terms 
and to control the growth of $\beta_{2,r}$ is to enforce, for $\ell <r$, 
some moderate uniform upper bound on $\|x_{\ell,j}-x_{\ell,0}\|$ only
depending on $\|s_{r,k}^L\|$, where iteration $(r,k)$ is the "ancestor" at 
level $r$ of iterations producing $x_{\ell,j}$ from $x_{\ell,0}$.  If this 
is acceptable (which, as we have mentioned, is not always beneficial in 
practice), the value of $\alpha_\ell$ in \req{alpha1-def} can be replaced by 
a smaller constant and the value of the second term in the definition of 
$\beta_{2,\ell+1}$ in \req{beta-rec12} then only grows moderately with $m_\ell$.
If this option is too restrictive, one may assume, or impose by a suitable 
termination criterion for lower-level iterations, that, for $\ell<r$, the 
gradients $g_{\ell,j}$ remain, on average, small enough in norm compared to
$\|g_{\ell,0}\|$ to ensure that
\[
 \EcondA{\ell+1,k}{\|s_{\ell,j}^L\|^2}
 \leq \kappa_g^2 \EcondA{\ell+1,k}{\|s_{\ell,0}^L\|^2}.
\]
Then this bound can be used in \req{xx}, 
avoiding the invocation of \req{bigs}.
As a consequence, $\alpha_\ell$ remains moderate
and the growth of $\beta_{2,\ell}$ in \req{beta-rec12} is also moderate.

In the deterministic case, Theorem~\ref{theorem:convergence} provides
a bound on the complexity of solving the bound-constrained problem
\req{problem} for which we consider the standard optimality measure
\[
\|\Xi_{r,k}\| \eqdef \|P_\calF(x_{r,k}-G_{r,k})-x_{r,k}\|
\]
(in the unconstrained case, $\|\Xi_{r,k}\|= \|G_{r,k}\|$).
Assumption AS.5 and AS.6 significantly simplify in this
context (that is when $g_{\ell,k} = G_{\ell,k}$ for all
$(\ell,k)$ and conditional expectations disappear. Indeed, AS.5
obviously always holds and,
unless the starting point is already first-order critical, $\calE$ may
always be made to happen in this context by suitably 
choosing $\varsigma$. It is also easy to enforce 
AS.6 at any given recursive iteration $(\ell+1,k)$ by replacing
$f_\ell(x)$ by 
\beqn{tau-corrected}
h_{\ell,k}(x) = (P_\ell^T G_{\ell+1,k}-g_{\ell,0})^T(x-x_{\ell,0}) + f_\ell(x),
\eeqn
in which case the left-hand side of \req{coherence-recursive} is
identically zero. This technique is the standard ``tau correction''
used in deterministic multigrid methods to ensure coherence of
first-order information between levels $\ell$ and $\ell-1$ (see
\cite{BrigHensMcCo00, kopanicakova2022use}, for instance). 

\lcor{theorem:convergence-deterministic}{
Suppose that AS.1--AS.4 hold, that the
\algname\ algorithm is applied to problem \req{problem} starting from
a non-critical $x_{r,0}$ with $\varsigma < \|d_{r,0}\|^2$, that
$g_{\ell,j} = G_{\ell,j}$ for all $(\ell,j)$ and that the
tau-correction is applied at each recursive iteration (i.e. \req{tau-corrected} is
always used). Then
\beqn{gradbound-deter}
\average_{j\in\iiz{k}} \|\Xi_{r,j}\|^2
\leq \frac{\kap{conv}}{k+1},
\eeqn
where the constant $\kap{conv}$ is computed as in
Theorem~\ref{theorem:convergence} using the values
\[
\kappa_\tau = \kappa_{Gg} = 0.
\]
}

\proof{
The choice 
$\varsigma < \|d_{r,0}\|^2$ 
(which is always possible) ensures
that $\calE$ occurs, the use of \req{tau-corrected} implies that
\req{coherence-recursive} holds with $\kappa_\tau = 0$, and
AS.5 always holds with $\kappa_{Gg}=0$  since
$g_{\ell,j}=G_{\ell,j}$. The result then follows from
Theorem~\ref{theorem:convergence}.
}

The more general stochastic case is more complicated because
Theorem~\ref{theorem:convergence} only covers the convergence of $\EA{\|d_{r,k}\|^2}$, which a first-order criticality measure of an
approximation (using $g_{r,k}$ instead of $G_{r,k}$) of problem \req{problem}.
However, as explained in \cite[Section~4]{BellGratMoriToin25}, the underlying stochastic distribution of the approximate gradient may not ensure that 
the rate of convergence of $\EA{\|d_{r,k}\|}$ given by \req{gradbound} 
automatically enforces a similar rate of convergence of 
$\EA{\|\Xi_{r,k}\|}$, the  relevant measure for problem \req{problem} itself.  
Fortunately, this can be alleviated if one is ready to make an assumption on the 
error of the gradient oracle.  A suitable assumption is obvioulsy to require that
\beqn{coherent}
\EA{\|\Xi_{r,k}\|}\leq \kap{opt} \EA{\|d_{r,k}\|},
\eeqn
for all $k\geq 0$ and some constant $\kap{opt}\in (0,1]$ (we then say that the gradient is "coherently distributed"). This condition is satisfied if the problem is 
unconstrained and the gradient oracle unbiased, because then, using Jensen's inequality and the convexity of the norm,
\[
\|\Xi_{r,k}\|
= \|G_{r,k}\|
= \|\EA{g_{r,k}}\|
\leq \EA{\|g_{r,k}\|}
= \EA{\|d_{r,k}\|}.
\]
When bounds are present, another suitable assumption is given by the next lemma.

\llem{lemma:bias}{
For each $\ell\in\iiz{r}$ and each $k\geq0$, we have that
\beqn{bias}
\EA{\|\Xi_{\ell,k}\| }
\leq \EA{\|d_{\ell,k}\|}+\EA{\|g_{\ell,k}-G_{\ell,k}\|}.
\eeqn
Moreover, if 
$\EcondA{\ell,k}{\|g_{\ell,k}-G_{\ell,k}\|}\leq \kap{bias} \,\EcondA{\ell,k}{\|d_{\ell,k}\|}$ for some $\kap{err}\geq 0$,
then
\beqn{unbiased}
\EA{\|\Xi_{\ell,k}\|} \leq (1+\kap{err})\EA{\|d_{\ell,k}\|}.
\eeqn
}

\proof{See \cite[Lemma~3.6]{BellGratMoriToin25}.}

We may now apply this result to derive a complexity bound for the
solution of problem \req{problem} from Theorem~\ref{theorem:convergence}.

\lthm{theorem:true-convergence}{
Suppose that AS.1--AS.6 hold and that the
\algname\ algorithm is applied to problem \req{problem}. 
Suppose also that, for all $k\geq 0$, either \req{coherent} holds or
\beqn{bias-condition}
\EcondA{r,k}{\|g_{r,k}-G_{r,k}\|}\leq \kap{err} \,\EcondA{r,k}{\|d_{r,k}\|^2}.
\eeqn
for some $\kap{err} \geq 0$. Then
\beqn{gradbound-true}
\average_{j\in\iiz{k}} \EA{\|\Xi_{r,j}\|^2}
\le \frac{\kap{conv}}{k+1},
\eeqn
with 
\beqn{kapconv-def-true}
\kap{conv}
=\bigfrac{1}{2}\varsigma\,\kappa_W^2(1+\kap{err})\,\left|W_{-1}\left(-\frac{1}{\kappa_W}\right)\right|^2
\leq \bigfrac{1}{2}\varsigma\,\kappa_W^2(1+\kap{err})\left|\log(\kappa_W) + \sqrt{2(\log(\kappa_w)-1)}\right|^2,
\eeqn
where $W_{-1}$ is the second branch of the Lambert function, 
$\Gamma_0 \eqdef f_r(x_0)-\flow$ and 
$\kappa_W$ is defined in \req{kappaW-def} with $\beta_{r,1}$ and $\beta_{2,r}$ 
constructed using the recurrence \req{beta-rec12}.
}

Note that the inequality \req{bias} in Lemma~\ref{lemma:bias} also indicates 
what happens if condition  \req{bias-condition} fails for large
$k$. In that case, the right-hand side of the inequality is no longer
dominated by its first term, which then only converges to the level 
of the gradient oracle's error.

We finally state  a complexity result in probability simply
derived from Theorem~\ref{theorem:true-convergence} by using Markov's inequality.

\lcor{prob-complexity}{Under the conditions of
  Theorem~\ref{theorem:true-convergence} and given $\delta\in
  (1-p_\calE,1)$ where  $p_\calE$ is the probability of occurence
  of the event $\calE$, one has that
\[
\Prob{\min_{j\in\iiz{k}} \|\Xi_{r,j}\| \leq
  \epsilon} \geq 1 -\delta
\tim{ for }
k \geq \frac{p_\calE \kap{conv}}{(p_\calE-(1-\delta))\epsilon^2}.
\]
}

\proof{See \cite[Corollary~2.6]{BellGratMoriToin25}.}

\noindent
Thus the \algname\ algorithm needs at most $\calO(\epsilon^{-2})$
iterations to ensure an $\epsilon$-approximate first-order critical
point of the bound-constrained problem \req{problem} with 
probability at least $1-\delta$, the constant being inversely 
proportional to $\delta$.
  
\section{Using additive Schwarz domain decompositions}\label{section:decomposition}

As we have pointed out in the introduction, the conditions we have 
imposed on $R_\ell$ and $P_\ell$ are extremely general. This section 
considers how this freedom can be exploited to handle domain-
decomposition algorithms.  That this is possible has already been 
noted, in the the more restrictive standard framework of 
deterministic methods using function values by \cite{gross2009unifying, GrosKrau11} and 
\cite{GratViceZhan19}. Following Section~5 of this last reference, 
we consider, instead of a purely hierarchical organization of the 
optimization problem \req{problem}, a description where subsets of 
variables describe the problem in some (possibly overlapping) 
"subdomains".  Our objective is then to separate the optimization 
between these subdomains as much as possible. 
 
More precisely, let $\{\calD_{{p}}  \}_{p=1}^M$ be a (possibly 
overlapping) \textit{covering} of $\ii{n}$, that is a 
collection of sets of indices such that
\[
\emptyset \not= \calD_p \subseteq \ii{n}
\tim{and}  \bigcup_{p=1}^M \calD_p = \ii{n}.
\]
When the problem is related to the discretization of a physical 
domain, it is often useful to assign to each $\calD_p$ the indices 
of the variables corresponding to discretization points in the 
physical subdomains, but we prefer a more general, purely algebraic definition. 
We say that a partition $\{\widehat{\calD}_p\}_{p=1}^M$ of $\ii{n}$ is a 
\textit{restriction\footnote{A restriction involves disjoint 
subsets of variables.} of the covering} $\{\calD_p\}_{p=1}^M$, if
$\widehat{\calD}_p\subseteq \calD_p$ for all $p\in \ii{M}$. We also 
assume the knowledge of linear mappings $\{R^{(p)}\}_{p=1}^M$ between 
$\Re^{n_r} = \Re^n$ and $\Re^{n_p}$ (where $n_p = |\calD_p|)$ and 
$\{P^{(p)}\}_{p=1}^M$ between $\Re^{n_p}$ and $\Re^{n}$.  Given the 
covering $\{\calD_p\}_{p=1}^M$ and one of its restrictions 
$\{\widehat{\calD}_p\}_{p=1}^M$, these operators may be defined in 
various ways.  If $\{e_j\}_{j=1}^n$ are the columns of the identity 
matrix on $\Re^n$, let the matrices $U_p$ and $\widehat{U}_p$ be 
given column-wise by 
\beqn{UphatUp-def}
U_p = \Big( \{e_j\}_{j\in \calD_p} \Big)
\tim{ and }
\widehat{U}_p = 
\left( \{e_j\}_{j\in\widehat{\calD}_p} \right),
\eeqn
(ordered by increasing value of $j$), and also define
\beqn{Wp-def}
W_p = \left[\sum_{q=0}^{M}U_q U_q^T\right]^{-1}U_p
\;= \, {\rm diag}(\theta_p)^{-1}U_p.
\eeqn 
Here, $\theta_p \in \Re^{n}$ is defined such that its $i$-th component equals to
the number of subdomains involving variable 
$i$ (see \cite[Lemma~E.1]{GratViceZhan19}).
Then the well-known additive Schwarz decompositions (see 
\cite{CaiSark99,FromSzyl02,Hols94}, for instance) 
can be defined by setting $P_p$, $R_p$ and $W_p$ 
as specified in Table~\ref{table:AS-defs}.
As can be seen in this table, and in contrast with the 
hierarchical case where the (Galerkin) choice
$R_\ell=P_\ell^T$ is often made, this relation does not hold 
for the considered decomposition methods, except for AS and 
RASH. 
We also note that the entries of $R_p$ and $P_p$ are non-negative.

\begin{table}[htbp]\small
\begin{center}
\begin{tabular}{llll}
Decomposition technique & Abbrev. & $P^{(p)}$ & $R^{(p)}$ \\
\hline
&&&\\*[0.1ex]
Additive Schwarz & AS & $U_p$ & $U_p^T$ \\
Restricted Additive Schwarz & RAS & $\widehat{U}_p$ & $U_p^T$ \\
Weighted Restricted Additive Schwarz & WRAS & $W_p$ & $U_p^T$ \\
Additive Schwarz (Harmonic) & ASH & $U_p$ & $\widehat{U}_p^T$ \\
Restricted Additive Schwarz (Harmonic) & RASH & $\widehat{U}_p$ & $\widehat{U}_p^T$\\
Weighted Additive Schwarz (Harmonic) & WASH & $U_p$ & $W_p^T$\\
\end{tabular}
\caption{\label{table:AS-defs}The prolongations and 
restrictions for the standard additive Schwarz 
domain decomposition techniques~\cite[Section~5.5.2]{GratViceZhan19}.}
\end{center}
\end{table}

Much as we did above for hierarchical models, we also associate a 
local model of the objective function 
$f^{(p)}:\Re^{n_p} \rightarrow \Re$ to each subdomain $\calD_p$.
We then consider the \textit{extended} space $\Re^{N}$ 
where $N= \sum_{p=1}^M n_p$ and define the operators
\beqn{RP-decomp}
R_1 = \left(\begin{array}{c} R^{(1)}\\ \vdots \\ R^{(M)}\end{array}\right)
\tim{ and }
P_1 = \left( P^{(1)}, \ldots, P^{(M)} \right)
\eeqn
from $\Re^n$ to $\Re^N$ and from $\Re^N$ to $\Re^n$, respectively. 
We also define the \textit{extended model}
by
\beqn{f0-def}
f_0( x^{(1)}, \ldots, x^{(M)}) = \sum_{p=1}^M f^{(p)}(x^{(p)}).
\eeqn 
This extended model is conceptually useful because our next step is 
to consider a two-levels hierarchical model of the type analyzed in 
the previous sections, where the "top" level ($\ell=1=r$) objective 
function is the original objective function $f_1=f$, while 
its lower level model $(\ell=0$) is the extended model 
$f_0$, the associated prolongation and 
restriction operators between $\Re^n$ and $\Re^N$ being given 
by $R_1$ and $P_1$ in \req{RP-decomp}. Note that all entries of 
$P_1$ are non-negative and $\|P_1\|_{\infty} = 1$.

We now investigate two interesting properties of this setting. 
The first (and most important) is that the minimization of $f_0$ is 
\textit{totally separable} in the variables $x^{(p)}$ and amounts 
to minimizing each $f_p(x^{(p)})$ independently.  In particular,
\textit{these minimizations may be conducted in parallel}.

The second is expressed in the following lemma.

\llem{lemma:simpler-lower-bounds}{
In all cases mentioned in Table~\ref{table:AS-defs}, the lower-level bounds given by
\beqn{simpler-low-bounds}
\hat{l}_0^{(p)} \eqdef R^{(p)} l
\tim{and}
\hat{u}_0^{(p)} \eqdef R^{(p)} u
\eeqn
coincide with those resulting from \req{lower-level-bounds}.  Moreover,
\beqn{even-simpler-low-bounds}
\hat{l}_0^{(p)} = l_{\calD_p}
\tim{ and }
\hat{u}_0^{(p)} = u_{\calD_p},
\eeqn
for AS, RAS and WRAS, and
\beqn{even-simpler-low-bounds-H}
\hat{l}_0^{(p)} = l_{\widehat{\calD}_p}
\tim{ and }
\hat{u}_0^{(p)} = u_{\widehat{\calD}_p},
\eeqn
for ASH and RASH.
}

\proof{
Let $p\in\ii{m}$ be fixed and observe first that, because $R^{(p)}$ is non-negative, 
$\hat{l}_0^{(p)}$ and $\hat{u}_0^{(p)}$ as defined by \req{simpler-low-bounds} are such that 
$\hat{l}_0^{(p)}\leq \hat{u}_0^{(p)}$. Observe also that $R_{j,i}^{(p)} >0$, if and only if 
$P_{i,j}^{(p)} >0$.

If the
$j$-th column of $P^{(p)}$ is zero, then the corresponding domain-dependent variable does not contribute to the prolongated step, i.e., $(P^{(p)}(x_{0,*}^{(p)}-x_{0,0}^{(p)}))_j = 0$ and $s_{1,k,j}=0$.
In particular, its bounds 
can be chosen as in \req{simpler-low-bounds}, \req{even-simpler-low-bounds} or 
\req{even-simpler-low-bounds-H}. Suppose now that column $j$ of $P^{(p)}$ is 
nonzero. Table~\ref{table:AS-defs} with \req{UphatUp-def} and \req{Wp-def} implies that it has 
only one nonzero entry, say in row $i$, so that \req{lower-level-bounds} gives that
\beqn{gen-bounds}
l_{0,j} = x_{0,0,j} + \frac{l_{1,i} - x_{1,k,i}}{\sigma_{1,i}}.
\eeqn
The relation between $R^{(p)}$ and $P^{(p)}$ and the inequality $P^{(p)}_{i,j}>0$  
also give that 
\beqn{partial}
[R^{(p)}l]_j = R^{(p)}_{j,i}l_{1,i}
\tim{ and }
x_{0,0,j}= [R^{(p)}x_{1,k}]_j = R^{(p)}_{j,i}x_{1,k,i}.
\eeqn
Suppose now that one of AS, RAS, WRAS, ASH or RASH is used.  Then $R^{(p)}_{j,i}=1$
and \req{partial} ensures that
\[
[R^{(p)}l]_j= l_{1,i} 
\tim{ and }
x_{0,0,j} = x_{1,k,i}.
\]
Substituting these equalities in \req{gen-bounds} and observing that 
$\sigma_{1,i}=1$ for all these variants then yields that
\beqn{l0j}
l_{0,j} 
= R^{(p)}_{j,i}x_{1,k,i} + R^{(p)}_{j,i}l_{1,i} - R^{(p)}_{j,i}x_{1,k,i} 
= R^{(p)}_{j,i}l_{1,i},
\eeqn
which proves the first part of \req{simpler-low-bounds} given that $l_1=l$. 
Moreover, \req{l0j} and the fact that $R^{(p)}_{j,i}=1$ when $j\in\calD_p$ for 
AS, RAS and WRAS imply that the first part of \req{even-simpler-low-bounds} 
holds. Similarly, \req{l0j} and the fact that $R^{(p)}_{j,i}=1$ when 
$j\in\widehat{\calD}_p$ for ASH and RASH ensure the first part of 
\req{even-simpler-low-bounds-H}.

If we now turn to WASH, we have from \req{Wp-def} that 
\[
R^{(p)}_{j,i}= \frac{1}{\theta_i} = \frac{1}{\sigma_{1,i}},
\]
where $\theta_i$ may now be larger than one, if variable $i$ occurs 
in more than one subdomain. Hence, using \req{partial},
\beqn{the-bad-one}
x_{0,0,j} = \frac{x_{1,k,i}}{{\sigma_{i,1}}}
\tim{ and }
\frac{l_{1,i}}{\sigma_{1,i}} = [R^{(p)}l]_j,
\eeqn
and \req{gen-bounds} ensures that the first part of \req{simpler-low-bounds} 
also holds for WASH\footnote{The second part of \req{the-bad-one} also implies 
that the first part of \req{even-simpler-low-bounds-H} fails for variables 
appearing in more than one subdomain.}.

Proving the result for the upper bounds uses the same argument.
}

\noindent
Capitalizing on these observations, we may now reformulate the 
\algname\ algorithm for our decomposition setting, as stated in 
Algorithm~\ref{top-algo-d} where ``recursion iterations'' now 
become ``decomposition iterations''.

\algo{top-algo-d}{\algnamed}{
\begin{description}
\item[Initialization:] The function $f_r=f_1=f$, bounds $l$ and $u$, $x_{1,0}$, such that 
$l_i < x_{1,0,i} < u_i$ for $i\in \ii{n}$, the decomposition and 
models $\{f^{(p)},P^{(p)},R^{(p)},m^{(p)}\}_{p=1}^m$, and the constants 
$m_1\geq0$, $\kappa_s,\kap{2nd} \geq 1$
and $\kap{1st},\kap{gs},\varsigma\in(0,1]$ are given.
\item[Top level optimization:]
Return {\sf \algnamed-r}$(1,x_{1,0},l,u,\varsigma^2, 0, \infty, m_1)$.
\end{description} 
}
\clearpage

\algo{rec-algo-d}{$x_{\ell,*} = $ {\sf \algnamed-r}
    ($\ell,x_{\ell,0},l_\ell,u_\ell,w_{\ell,-1}, \theta_{1,\ell},\theta_{2,\ell}, m_\ell$)}{
\begin{description}
\item[Step 0: Initialization:] Set $k=0$. 
\item[Step 1: Start iteration:]
  Identical to Step~1 of the \algname-r algorithm, where iterations at level $\ell=1$ are either "Taylor" or "decomposition".

\item[Step 2: Taylor iteration: ] 
Identical to Step~2 of the \algname-r algorithm.

\item[Step 3: Decomposition iteration: ]  Set
  \beqn{thetas-def-d}
  \theta_{1,1} = \kap{1st} \,|d_{1,k}^T\Delta_{1,k}|,
  \ms\ms
  \theta_{2,1}= \kap{2nd} \|s_{1,k}^L\|.
  \eeqn  
   Then, for each $p\in\ii{m}$, compute
  \beqn{kappadD-def-d}
  [ x_{0,0}^{(p)}, w_0^{(p)}, l_0^{(p)}, u_0^{(p)}] 
  =  R^{(p)} [ x_{1,k}, w_{1,k}, l, u ],
  \eeqn
  and set
  \beqn{parallel}
  \hspace*{-5mm}
  s_{0,k}^{(p)} = P^{(p)} \Big[ \mbox{\sf \algnamed-r}(0,x_{0,0}^{(p)},l_0^{(p)},u_0^{(p)},w_0^{(p)},\theta_{1,1}, \theta_{2,1}, m^{(p)}) -x_{0,0}^{(p)}\Big].
  \eeqn
  Then set
  \beqn{synchro}
  s_{1,k} = \sum_{p=1}^m s_{0,k}^{(p)}.
  \eeqn
\item[Step 4: Loop:] 
  Identical to Step~4 of the \algname-r algorithm.
\end{description} 
}

\noindent
A mentioned above, the most attractive feature of the \algnamed\ 
algorithm is the possibility to compute the recursive 
domain-dependent $m$ minimizations of Step~3 (i.e.\ \req{kappadD-def-d} 
and \req{parallel}) in parallel, 
the synchronization step \req{synchro} being particularly simple. 
But, as for \algname, the choice between a Taylor and a 
decomposition iteration remains open at top level iterations, 
making it possible to alternate iterations of both types.
Note also that, because of the second part of 
Lemma~\ref{lemma:simpler-lower-bounds}, \req{kappadD-def-d} is equivalent to 
(and thus may be replaced by)
\[
[ x_{0,0}^{(p)}, w_0^{(p)}] =  R^{(p)} [ x_{1,k}, w_{1,k}],
\ms
l_0^{(p)} = l_{\calD_p}
\tim{ and }
u_0^{(p)} = u_{\calD_p}
\]
when one of AS, RAS or WRAS is used, or to
\[
[ x_{0,0}^{(p)}, w_0^{(p)}] =  R^{(p)} [ x_{1,k}, w_{1,k}],
\ms
l_0^{(p)} = l_{\widehat{\calD}_p}
\tim{ and }
u_0^{(p)} = u_{\widehat{\calD}p}
\]
for ASH or RASH, vindicating our comment in the introduction on the 
ease of implementing the constraint restriction rules for the 
domain-decomposition methods.

Observe now that, if \req{nogain} holds in Step~1 of the 
minimization in the $p$-th subdomain (in Step~3 of \algnamed), 
then $s_{0,k}^{(p)} =0$. Thus the loop on $p\in\ii{m}$ is Step~3 
and the sum in \req{synchro} are in effect restricted to the 
indices belonging to
\beqn{Ik-def}
\calI_k = \left\{ p \in\ii{m} \mid \left|\left(d_{0,0}^{(p)}\right)^T\Delta_{0,0}^{p}\right| \geq \theta_{1,0} = \frac{\kap{1st}}{m}\left|d_{1,k}^T\Delta_{1,k}\right|\right\},
\eeqn
which corresponds to choosing $f_0$ as
\[
f_0( x^{(1)}, \ldots, x^{(m)}) = \sum_{p\in\calI_k} f^{(p)}(x^{(p)})
\]
instead of \req{f0-def} (using the freedom to make the lower-level model 
depend on the upper level iteration).
Then, if $q = |\calI_k|$, $\calI_k = \{p_1,\ldots,p_q\}$, 
\[
d_{0,0}^T = \left(\left(d_{0,0}^{(p_1)}\right)^T, 
\ldots, \left(d_{0,0}^{(p_q)}\right)^T\right)
\tim{ and }
\Delta_{0,0}^T = \left(\left(\Delta_{0,0}^{(p_1)}\right)^T, 
\ldots, \left(\Delta_{0,0}^{(p_q)}\right)^T\right),
\]
we have from \req{thetas-def-d} that 
\[
|d_{0,0}^T \Delta_{0,0}| 
= \sum_{p\in\calI_k} |\left(d_{0,0}^{(p)}\right)^T\Delta_{0,0}^{p}|
\geq \kap{1st}\,|d_{1,k}^T\Delta_{1,k}|,
\]
and thus the linear decrease for the model $f_0$ is at least 
a fraction\footnote{Because the requested decrease is now on a single subdomain, it makes sense to choose $\kap{1st}$ smaller than would would be requested for the full domain.} of $|d_{1,k}^T\Delta_{1,k}|$, as requested by \req{nogain}. 
We also have that
\[
\theta_{2,0}
= \sum_{p\in \calI_k}\theta_{2,0}^{(p)}
\leq m\,\kap{2nd}\|s_{1,k}^L\|, 
\]
and \req{thetas-def} therefore holds with $\kap{2nd}$ replaced by $m\,\kap{2nd}$
We conclude from this discussion that we may
apply Theorems~\ref{theorem:true-convergence} and 
Corollaries~\ref{theorem:convergence-deterministic} and 
\ref{prob-complexity} to the \algnamed\ algorithm, provided we review our assumptions in the domain decomposition context.
Fortunately, this is very easy, 
as it is enough to assume that AS.2 to AS.6 now hold for 
$\ell=1$ and each $p\in\ii{m}$ instead of for each 
$\ell\in\iiz{r}$. Thus these assumptions hold for 
level $1$ and for the extended problems of level $0$.
As a consequence, we also obtain that the \algnamed\ algorithm 
converges to a first-order critical point with a rate prescribed 
by these results. In particular, Corollary~\ref{prob-complexity} 
ensures that, with high probability, the \algnamed\ decomposition 
algorithm has an $\calO(\epsilon^{-2})$ evaluation complexity to 
achieve an $\epsilon$-approximate critical point.

\section{Numerical results}
\label{section:numerics}

In this section, we numerically study the convergence behavior of the proposed \algname\ and \algnamed\ algorithms.

\subsection{Benchmark problems}
We utilize the following five benchmark problems, encompassing both PDEs and deep neural network (DNNs) applications.\\

\begin{figure}[htb]
\includegraphics[scale=0.05]{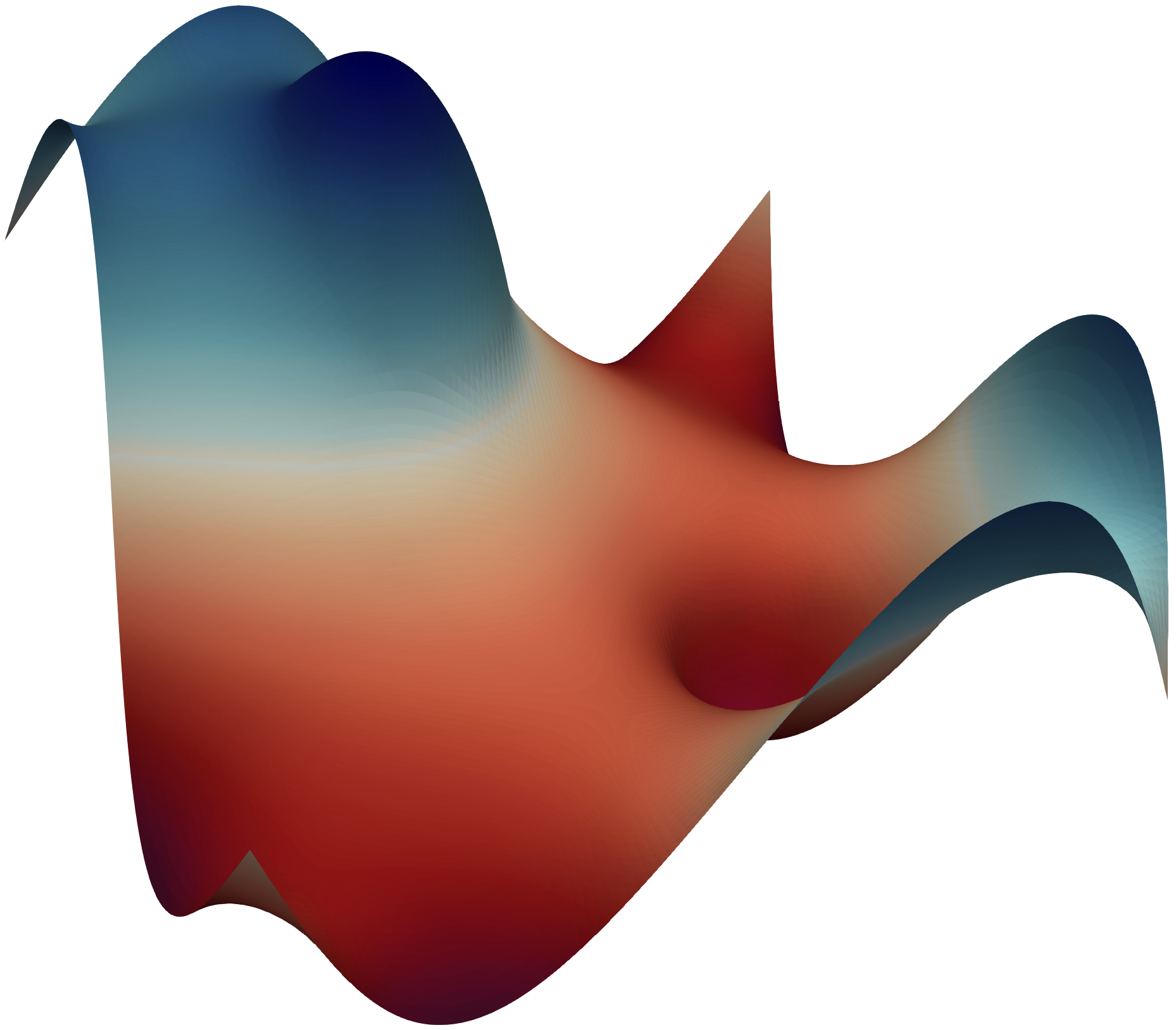}
\includegraphics[scale=0.1]{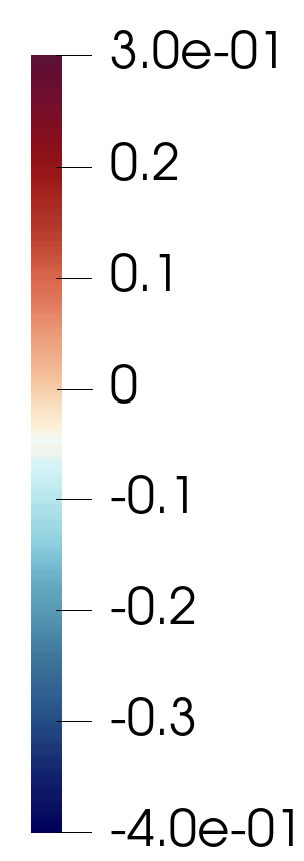}
\hfill
\includegraphics[scale=0.05]{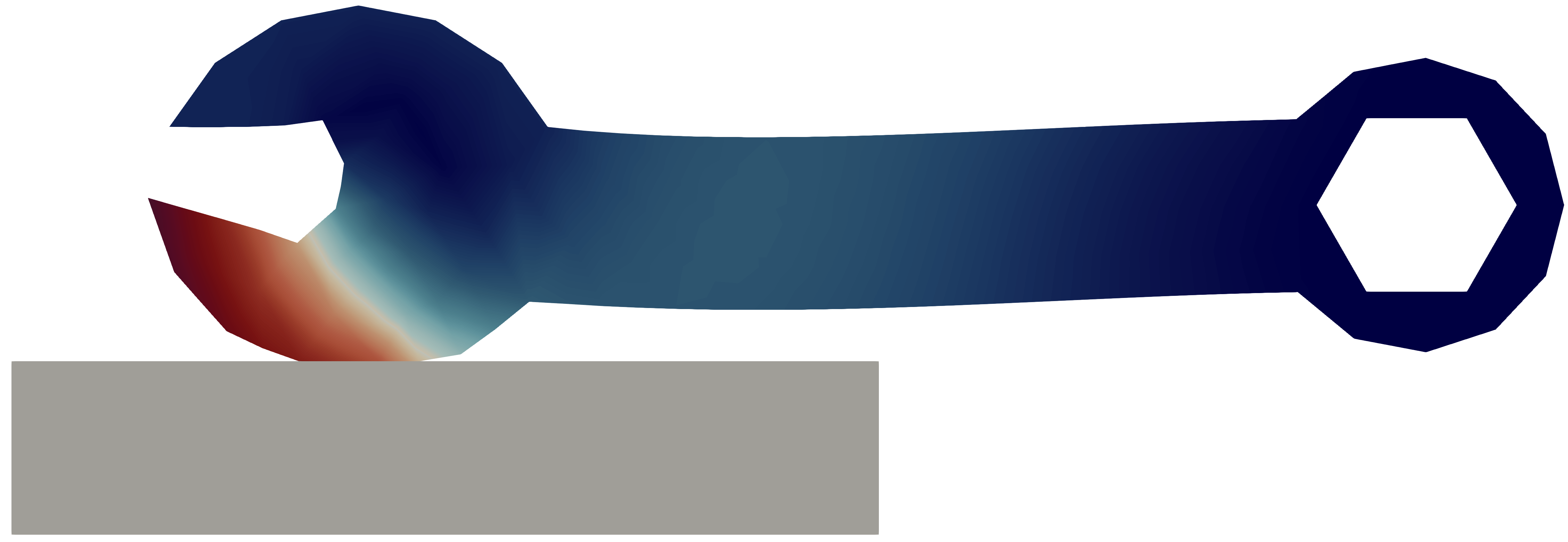}
\includegraphics[scale=0.05]{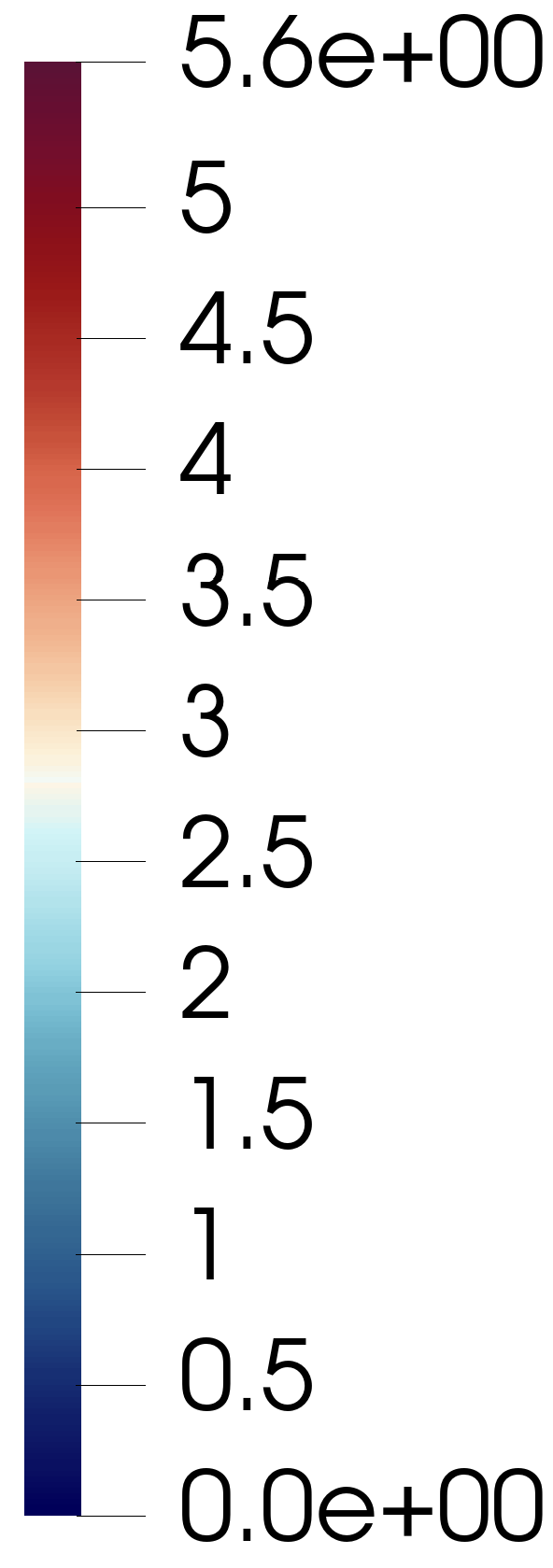}
\caption{Simulation results (solution $z$) for \textit{Minsurf} and \textit{NeoHook} examples.}
\label{fig:sim_results}
\end{figure}

\noindent\textbf{\textit{1.~Membrane}:}
Following~\cite{domoradova2007projector}, let ${\Omega:=(0,1)^2}$ be a computational domain with boundary~${\Gamma=\partial \Omega}$, decomposed into three parts:~${\Gamma_l = \{0\} \times (0,1)}$,~${\Gamma_r = \{1\} \times (0,1)}$, and~${\Gamma_f = (0,1) \times \{0,1\}}$. 
The minimization problem is given as
\begin{equation}
\begin{aligned}
& \underset{z \in \pazocal{Z}}{\text{min}} \ \ f(z) := \frac{1}{2} \int_{\Omega} \| \nabla z(x) \|^2 \ d z + \int_{\Omega} z(x) \ d z, \\
& \text{subject to} \ \ l(x) \leq z, \quad \text{on~$\Gamma_r$},
\end{aligned}
\label{eq:min_membrane}
\end{equation}
where~$\pazocal{Z}:=\{z \in H^1(\Omega) \ | \ z = 0 \ \text{on} \ \Gamma_l \}$.
The lower bound~$l$ is defined on~$\Gamma_r$, by the upper part of the circle with the radius,~$R=1$, and the center,~$C=(1,-0.5,-1.3)$, i.e., 
\begin{align*}
l(x) = 
\begin{cases} (-2.6 + \sqrt{2.6^2 - 4((x_2 - 0.5)^2 - 1.0 + 1.3^2) })/2,  & \quad \text{if} \  x_1=1, \\
-\infty, &  \quad \text{otherwise},
\end{cases}
\end{align*}
where the symbols~$x_1, x_2$ denote spatial coordinates.
The problem is discretized using a uniform mesh with 120 x 120 Lagrange finite elements ($\mathbb{Q}_1$).

\noindent\textbf{\textit{2.~MinSurf}:}
Let us consider the minimal surface problem~\cite{kovcvara2016first}, given as
\begin{equation}
\begin{aligned}
& \underset{z \in \pazocal{Z}}{\text{min}} \ \  f(z) := \int_{\Omega} 1 +  \| \nabla z(x) \|^2  \ d z, \\
& \text{subject to} \ \  l(x) \leq z \leq u(x), \quad  \text{a.e. in~$\Omega$},
\end{aligned}
\label{eq:min_sruf}
\end{equation}
where ${\Omega:=(0,1)^2}$ and the variable bounds are defined as
\begin{equation*}
\begin{aligned}
l(x) 	&= 0.25 - 8(x_1 - 0.70)^2 - 8 (x_2 - 0.70)^2,    \\
u(x) 	&= -(0.4 - 8 (x_1 - 0.3)^2 - 8 (x_2 - 0.3)^2).
\end{aligned}
\end{equation*}
The boundary conditions are prescribed such that~$\pazocal{Z}:=\{z \in H^1(\Omega) \ | \ z = g \ \text{on} \ \partial \Omega \}$, where $g$ is specified as
\begin{align}
g(x) = 
\begin{cases}
- 0.3 \sin(2 \pi x_2),  & x_1=0, \  0 \leq x_2 \leq 1, \\
+ 0.3 \sin(2 \pi x_2),  & x_1=1, \ 0 \leq x_2 \leq 1, \\
- 0.3 \sin(2 \pi x_1),  & x_2=0, \ 0 \leq x_1 \leq 1, \\
+ 0.3 \sin(2 \pi x_1),  & x_2=1, \ 0 \leq x_1 \leq 1. 
\end{cases}
\end{align}
The problem is discretized using a uniform mesh with $120\times 120$ Lagrange finite elements ($\mathbb{P}_1$).

\noindent\textbf{\textit{3.~NeoHook}:}
Next, we consider a finite strain deformation of a wrench made of rubbery material. 
Let~$\Omega$ be a computational domain, which represents a wrench of length $80$ mm, width of $24$ mm and thickness of $2$ mm, see also Figure~\ref{fig:sim_results} on the right. 
We employ the Neo-Hookean material model, and seek the displacement~${z}$ by solving the following non-convex minimization problem: 
\begin{equation}
\begin{aligned}
& \underset{{z} \in \pazocal{Z}}{\text{min}} \ \  f(z) := \frac{1}{2} \int_{\Omega} \frac{\mu}{2} (I_C - 3) - \mu (\ln(J)) + \frac{\lambda}{2} (\ln(J))^2 \ d z  - \int_{\Gamma_{top}} q_2 z_2 \ d s, \\
& \text{subject to} \ \  {g(x)} \geq {0}, \quad  \text{a.e. in~$\Omega$},
\end{aligned}
\label{eq:neohook}
\end{equation}
where~$\pazocal{Z}:=\{{z} \in [H^1(\Omega)]^3 \ | \ z = 0 \ \text{on} \ \Gamma_{in} \}$.
The symbol $\Gamma_{in}$ describes the part of the boundary corresponding to the inner part of the jaw, placed on the right side.
Moreover, $J=\text{det}({F})$ denotes the determinant of the deformation gradient~${F} = {I} + \nabla {z}$. 
The first invariant of the right Cauchy-Green tensor is computed as $I_C := \text{trace}({C})$, where ${C} = {F}^T {F}$.
For our experiment, the Lamé parameters $\mu = \frac{E}{2(1+\nu)}$ and $\lambda=\frac{E \nu}{(1+\nu)(1-2\nu)}$ are obtained by setting the value of Young's modulus $E=10$ and Poisson's ratio $\nu = 0.3$. 
The last term in~\eqref{eq:neohook} corresponds to Neumann boundary conditions applied on~${\Gamma_{top}}$, which corresponds to the top of the left wrench's jaw. 
The traction ${q_2:=-2.5}$ is prescribed in $x_2$-direction and $z_2$ is the second component of the displacement vector~$z$.

The constraint is setup such that the body~$\Omega$ does not penetrate the surface of the obstacle~$\Gamma_{obs}:=\{ {x} \in \Re^3 \ | \  x_2 = -9.0 \}$. 
Thus, the gap function is defined as ${g(x)}:= \langle \psi(x) - x_{obs}, {n} \rangle$, where $\psi(x)= {x} + {z}(x)$ specifies the deformed configuration. 
Here, the unit outward normal vector is defined as~${n} = (0, -1, 0)^T$, and $ x_{obs}  = (x_1, -9.0, x_3)^T$. 
The problem is descritized using an unstructured mesh with Lagrange finite elements ($\mathbb{P}_1$), giving rise to a problem with~$45,016$ dofs. 

\noindent\textbf{\textit{4.~IndPines}:}
This example focuses on the soil segmentation using hyperspectral images provided by the Indian Pines dataset~\cite{baumgardner2015220}. 
Let~${\pazocal{D} = \{ (y_s, c_s) \}_{s=1}^{n_s}}$ be a dataset of labeled data, where~${y_s \in \Re^{n_{in}}}$ represents input features and~${c_s \in \Re^{n_{out}}}$ denotes a desirable target. 
Following~\cite{queiruga2020continuous,chang2017multi}, we formulate the supervised learning problem as the following continuous optimal control problem~\cite{haber2017stable}:
\begin{align}
 &\underset{Q, q, \theta, W_T, b_T}{\text{min}}
\ \ f(Q, q, \theta, W_T, b_T) :=  \frac{1}{n_{s}} \sum_{s=1}^{n_{s}} g(\pazocal{P}(W_T {q}_s(T) + b_T), c_s) 
+  \int\limits_{0}^{T} \pazocal{R}({\theta}(t) ) \ dt 
 +  \pazocal{S}({W}_T, b_T) ,
 \nonumber \\  
 & \text{subject to} \quad \partial_t {q}_s(t) = \pazocal{F}({q}_s(t), {\theta}(t)), \qquad \forall t \in (0,T], \label{eq:cts_problem} \\
& \quad \quad \quad \quad \quad \quad {q}_s(0) = Qy_s, \nonumber
\end{align}
where $q$ denotes time-dependent states from $\Re$ into $\Re^{n_{fp}}$ and ${\theta}$ denotes the time-dependent control parameters from  $\Re$ into  $\Re^{n_{c}}$.
The constraint in~\eqref{eq:cts_problem} continuously transforms an input feature~${y_s}$ into final state~${q}_s(T)$, defined at the time $T$.
This is achieved by firstly mapping the input $y_s$ into the dimension of the dynamical system as~${q}_s(0) = Qy_s$, where~${{Q} \in \Re^{{n_{fp}} \times n_{in}}}$. 
Secondly, the nonlinear transformation of the features is performed using a residual block~$\pazocal{F}({q}_s(t), {\theta}(t)):=\sigma(W(t) q_s(t) + b(t))$, where ${\theta(t) = (\text{flat}(W(t)), \text{flat}(b(t)))}$,
$\sigma$ is an \textit{ReLu} activation function from~$\Re^{{n_{fp}}}$ into~$\Re^{{n_{fp}}}$, ${b(t) \in \Re^{n_{fp}}}$ is the bias and ${W(t) \in \Re^{n_{fp} \times n_{fp}}}$ is a dense matrix of weights.
We set the layer width $n_{fp} $ to be $50$. 

We employ the softmax hypothesis function ($\pazocal{P}$ from~$\Re^{n_{out}}$ into~$\Re^{n_{out}}$) together with the cross-entropy loss function, defined as~${{g}(\hat{c}_s, c_s) := c_s^T \log (\hat{c}_s)}$, where~$\hat{c}_s := \pazocal{P}(W_T {q}_s(T) + b_T) \in \Re^{n_{out}}$. 
The linear operators $W_T \in \Re^{n_{out} \times n_{fp}}$, and $b_T \in \Re^{n_{fp}}$ are used to perform an affine transformation of the extracted features~${q}_s(T)$.
Furthermore, we utilize a Tikhonov regularization,  i.e.,~$\pazocal{S}(W_T, b_T):= \frac{\beta_1}{2} \| W_T \|^2_F + \frac{\beta_1}{2} \| b_T \|^2$, where $\| \cdot \|_F$ denotes the Frobenius norm. 
For the time-dependent controls, we use $\pazocal{R}(\theta(t)) :=  \frac{\beta_1}{2} \| \theta(t) \|^2  + \frac{\beta_2}{2} \| \partial_t \theta(t) \|^2$, where the second term ensures that the parameters vary smoothly in time, see~\cite{haber2017stable} for details.  

To solve the problem~\eqref{eq:cts_problem} numerically, we use the forward Euler discretization with equidistant grid ${0 = \tau_0 < \cdots < \tau_{K} = T}$, consisting of $K+1$ points. 
The states and controls are then approximated at a given time~$\tau_k$ as~${{q}_k \approx {q}(\tau_k)}$, and ${{\theta}_k \approx {\theta}(\tau_k)}$, respectively.
Note, each~$\theta_k$ and~${q}_k$ now corresponds to parameters and states associated with the $k$-th layer of the ResNet DNN. 
The numerical stability is ensured by employing a sufficiently small time-step~$\Delta_t = T/(K)$. 
In particular, we set $K, T, \beta_1, \beta_2$  to $K=17$, $T=3$ and $\beta_1=\beta_2= 10^{-3}$.

\begin{figure}
\begin{minipage}{0.45\linewidth}
\scalebox{.45}{
\includegraphics{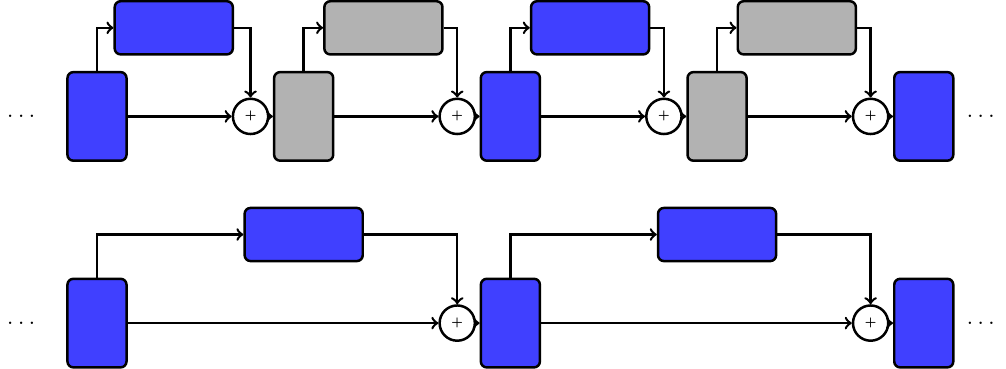}
}
\end{minipage}    
\hfill
\begin{minipage}{0.48\linewidth}    
\scalebox{0.48}{	
\includegraphics{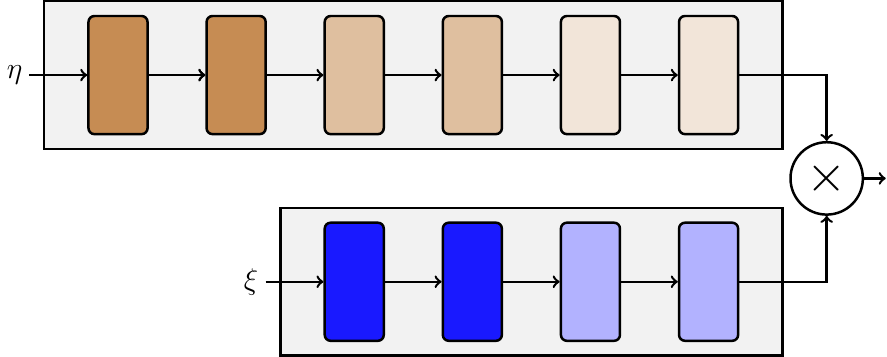}
}
\end{minipage}
\caption{Left: An example of a multilevel hierarchy of ResNets, with the upper level depicted on the top of the figure.
Right: An example of decomposing the parameters of DeepONet with branch (brown) and trunk (blue) sub-networks into five subdomains. Each subdomain is depicted in a different color.}    
\label{fig:hierarchies}
\end{figure}

\noindent\textbf{5.~\textit{Aniso}:}
This example considers operator learning using DeepONet~\cite{lu2021learning} in order to approximate a solution of the following parametric anisotropic Poisson's equation: 
	      \begin{equation}
		      \begin{aligned}
			      - \nabla \cdot (K({\eta}) \ \nabla z(x) ) &= f(x),  \quad \quad  &  & \forall \ x \in \Omega, \\
			      z                                                           & = 0, \      &  & \forall \ x \in  \partial \Omega, \\
		      \end{aligned}
        \label{eq:aniso_laplace}
	      \end{equation}
where~$\Omega := (-1,1)^2$  is the computational domain, with boundary $\partial \Omega$.
The symbol~$z$ denotes  the solution and $f(x):=1$.
The diffusion coefficient~$K(\eta)$ is the following anisotropic tensor:
\begin{align*}
K(\eta) = 
\begin{pmatrix}
\cos(\alpha) & -\sin(\alpha) 	\\
\sin(\alpha) & \ \ \ \cos(\alpha) 	
\end{pmatrix}
\begin{pmatrix}
1 & 0 	\\
0 & \beta	
\end{pmatrix}
\begin{pmatrix}
\ \ \ \cos(\alpha) & \sin(\alpha) 	\\
-\sin(\alpha) & \cos(\alpha) 	
\end{pmatrix},
\end{align*}
parametrized using~${\eta}:=[\alpha, \beta]$. 
We sample the anisotropic strength~${\beta \in [10^{-6},1]}$ according to ${\log_{10} (\sfrac{1}{\beta}) \sim \pazocal{U}[0, 6]}$.
The angle of the anisotropic direction, denoted by the parameter~${\alpha \in (0,\pi)}$, is sampled as~$\alpha \sim \pazocal{U}[0, \pi]$.
To construct the dataset, for each value of ${\eta}$, we discretize~\eqref{eq:aniso_laplace} using a uniform mesh with $34 \times 34$ finite elements.
The dataset~$\pazocal{D} = \{ (\eta_s, \xi_s, y_s) \}_{s=1}^{n_s}$ is composed of $4{,}250$ training samples and $750$ testing samples.
Here, $\xi_s \in \mathbb{R}^{n_c \times 2}$ denotes the set of coordinates at which the high-fidelity finite element solution~$y_s \in \mathbb{R}^{n_c}$ is evaluated.

The training process is defined as the minimization between the DeepONet output and  the target solution, i.e., 
\begin{align}
\min_{z \in \mathbb{R}^n}  f(z) :=  \frac{1}{n_s n_c} \sum_{s=1}^{n_s} \sum_{c=1}^{n_c} \big( \langle B(\eta_{s}; z), T(\xi_{s,c}; z) \rangle - y_{s,c} \big)^2,
\label{eq:DON_min}
\end{align}
where $z$ denotes collectively parameters of the DeepONet, consisting of branch (B) and trunk (T) sub-networks. 
Their outputs $B(\eta_{s}; z) \in \mathbb{R}^{Q}$ and $T(\xi_{s,c}; z) \in \mathbb{R}^{Q}$ are linearly combined by means of the inner product to approximate the solution of~\eqref{eq:aniso_laplace} at the location~$\xi_{s,c}$ for a given set of parameters~$\eta_s$.
Both sub-networks have dense feedforward architecture with a width of 128, \textit{ReLu} activation functions, and an output dimension of $Q = 126$.
The branch consists of 6 layers, while trunk network has 4 layers.

\subsection{Implementation and algorithmic setup}
For the FEM examples, the benchmark problems and \algname \  are implemented using Firedrake~\cite{rathgeber2016firedrake}, while \algnamed \ is implemented using MATLAB. 
For the DNNs and their training procedures using \algname\ and \algnamed \ algorithms, we utilize PyTorch~\cite{imambi2021pytorch}.

Our \algname\ is implemented in the form of a V-cycle.
For \algname, we employ three Taylor iterations on all levels for both pre-and post-smoothing.
On the lowest level, i.e., $ \ell=0$, we perform five Taylor steps.
For \algnamed, we perform ten decomposition iterations before each Taylor iteration, for all examples except \textit{Aniso}, where we employ three Taylor iterations and 25 decomposition iterations. For FEM examples, we computed the product $(s_{\ell,k}^L)^TB_{\ell,k}s_{\ell,k}^L$ appearing in \req{gamma-def} using gradient finite-differences in complex arithmetic~\cite{abreu2018accuracy}, a technique which produces cheap machine-precision-accurate approximations of $(s_{\ell,k}^L)^TH_{\ell,k}s_{\ell,k}^L$ without any evaluation of the Hessian $H_{\ell,k}$.
For DNNs examples, we consider purely first-order variants of \algname\ and \algnamed, i.e., $B_{\ell,k} = 0$ for all $(\ell,k)$.
Unless specified otherwise, all evaluations of gradients are exact, i.e., without noise.
Parameters~$\varsigma$, $\kappa_{\text{1st}}$ and $\kappa_{\text{2nd}}$ are set to $\varsigma = 0.01$, 
$\kappa_{\text{1st}}=0.95$, $\kappa_{\text{2nd}}=10$ for all numerical examples.
The learning rate is set to $10^{-2}$ for \textit{IndPines} (on all levels) and $5 \times 10^{-2}$ for \textit{Aniso} (on all subdomains), as suggested by both single-level and multilevel/domain-decomposition hyperparameters tuning.

The hierarchy of objective functions~$\{f_{\ell}\}_{\ell=0}^r$ required by the \algname\ is constructed for FEM examples by discretizing the original problem using meshes coarsened by a factor of two.
Similarly, for the \textit{IndPines} example, we obtain a multilevel hierarchy by coarsening in time by a factor of two, i.e., each~$f_{\ell}$ is then associated with a network of different depths.
For all examples, the transfer operators~$\{P_{\ell}\}_{\ell=0}^{r-1}$ are constructed using piecewise linear interpolation in space or time; see~\cite{Kopanicakova_2020c, kopanicakova2022globally} for details regarding the assembly of transfer operators for ResNets.
The restriction operators~$\{R_{\ell}\}_{\ell=0}^{r-1}$ are obtained as $R_{\ell} = 2^{-d} P_{\ell}^T$, where $d$ denotes the dimension of the problem at hand.
Moreover, for FEM examples, we enforce the first-order consistency relation~\eqref{tau-corrected} on all levels, while no modification to~$\{f_{\ell}\}_{\ell=0}^r$ is applied for the DNNs examples. 

In the case of \algnamed \ algorithm, the function  $f_0$ is as specified in~\eqref{f0-def}.
For FEM examples, each $f^{(p)}$ is created by restricting $f_r$ to a subdomain determined by the METIS partitioner~\cite{karypis1997metis}. 
Transfer operators are assembled using the WRAS technique (see Table~\ref{table:AS-defs}) for the different overlap sizes in our tests.  
For the \textit{Aniso} example, the decomposition is obtained by first decomposing the parameters of branch and trunk networks into separate subdomains.
The parameters of both sub-networks are further decomposed in a layer-wise manner; see~\cite {kopanivcakova2024enhancing, lee2024two} for details regarding the construction of the transfer operators.
Figure~\ref{fig:hierarchies} illustrates multilevel and domain-decomposition hierarchies for our DNNs examples.  

For all FEM examples, the initial guess is prescribed as a vector of zeros, which satisfies the boundary conditions and does not violate the bound constraints.
For the DNNs examples, the network parameters are initialized randomly using Xavier initialization \cite{GlorBeng10}.
Therefore, all reported results for these set of examples are averaged over ten independent runs.

The optimization process for all FEM examples is terminated as soon as
$\|\Xi_{r,k}\| < 10^{-7}$ or $\|\Xi_{r,k}\| / \|\Xi_{r,0}\| < 10^{-9}$.
For the \textit{IndPines} example, the algorithms terminate as soon as $\text{acc}_{\text{train}} > 0.98$ or $\text{acc}_{\text{val}} > 0.98$. 
Termination also occurs as soon as $\sum_{i=1}^{15}
(\text{acc}_{\text{train}})_e - (\text{acc}_{\text{train}})_{e-i} <
0.001$ or $\sum_{i=1}^{15} (\text{acc}_{\text{val}})_e -
(\text{acc}_{\text{val}})_{e-i} < 0.001$, where
$(\text{acc}_{\text{train/val}})_e$ is defined as the ratio between
the number of correctly classified samples from the train/validation
dataset and the total number of samples in the train/validation
dataset for a given epoch~$e$.
For the \textit{Aniso} example, the algorithms terminate if the number of epoch exceeds $1,000,000$ or as soon as  $\sum_{i=1}^{100}
(f_{\text{train}})_e - (f_{\text{train}})_{e-i} <
0.0001$ or $\sum_{i=1}^{100} (f_{\text{val}})_e -
(f_{\text{val}})_{e-i} < 0.0001$, where $f$ denotes the loss function, as defined in~\eqref{eq:DON_min}. 

\subsection{Numerical performance of \algname \ }
We first investigate the numerical performance of multilevel \algname\ compared with its single level variant \algnameSL.
To assess the performance of the proposed methods fairly, we report not only the number of V-cycles, but also the dominant computational cost associated with gradient evaluations. 
Let $\pazocal{C}_r$ be a computational cost associated with an evaluation of the gradient\footnote{We recall that the Hessian is never evaluated, even when the second-order stepsize $\gamma_{\ell,k}$ is computed.} on the uppermost level,
the total computational cost~$\pazocal{C}_{\text{ML}}$ for \algname\  is evaluated as follows:
\begin{align}
 \pazocal{C}_{\text{ML}} =\sum_{\ell=0}^r \frac{n_{\ell}}{n_{r}} \sharp_\ell \, \pazocal{C}_r,
 \label{eq:cost_ML}
\end{align}
where~$\sharp_\ell$ describes a number of evaluations performed on a level~$\ell$, while the scaling factor~$\frac{n_{\ell}}{n_{r}}$ accounts for the difference between the cost associated with level~$\ell$ and the uppermost level~$r$. 
Note that for the FEM examples, this estimate accounts for the cost associated with the gradient evaluation required by the finite-difference method used to compute the product~$(s_{\ell,k}^L)^T B_{\ell,k} s_{\ell,k}^L$.
Moreover, we point out that in the case of DNN applications, additional scaling must be applied to account for the proportion of dataset samples used during evaluation.

We begin our numerical investigation by comparing the cost $\pazocal{C}_{\text{ML}}$ of \algname\ with \algnameSL. 
Table~\ref{tab:ML} presents the numerical results obtained for four benchmark problems, where increasing the number of refinement levels also corresponds to use more levels by \algname.
As our results indicate, the cost of \algname\ is significantly lower than that of \algnameSL \ for all benchmark problems. 
Notably, the observed speedup achieved by \algname\ increases with the number of levels, highlighting the practical benefits of the proposed algorithm for large-scale problems.
Moreover, for the \textit{IndPines} examples, \algname\ produces DNN models with accuracy that is comparable to, or even exceeds, that of \algnameSL. 
Based on our experience, this improved accuracy can be attributed to the fact that the use of hierarchical problem decompositions tends to have a regularization effect~\cite{GratKopaToin23}.

\begin{table}[]
\centering
\begin{tabular}{|r||r|r||r|r|r|r|r|r}
\hline
  \multirow{2}{*}{ Example} &   \multicolumn{2}{c||}{\multirow{2}{*}{Method}} & \multicolumn{4}{c|}{Refinement/\algname \  levels} \\ \cline{4-7}
  & \multicolumn{1}{c}{ }    & \multicolumn{1}{c||}{ } & \textbf{2} 	& \textbf{3} 	& \textbf{4} 	& \textbf{5}  		\\ \hline \hline
  \multirow{2}{*}{ Membrane} & \multicolumn{2}{r||}{\algnameSL}        & 2,704       	& 9,968      	&       36,054     		&     58,880    	    \\ \cline{2-7} 
& \multicolumn{2}{r||}{{\algname }}        						& 560        	& 588        	&        944    		&   2,102    		   \\ \hline \hline

  \multirow{2}{*}{ MinSurf} & \multicolumn{2}{r||}{{\algnameSL}}        		&	2,564  	&    10,400   	&     41,136     &  103,458       \\ \cline{2-7} 
&\multicolumn{2}{r||}{ {\algname } }        							& 684        	& 1,142        	&        2,448    		&   5,224    	   \\ \hline \hline

  \multirow{2}{*}{NeoHook} & \multicolumn{2}{r||}{\algnameSL}        	&	  67,704	&    229,512   	&      632,718  &     $> 750,000$    \\ \cline{2-7} 
& \multicolumn{2}{r||}{{\algname }}         						&    3,272	&   4,272      	&   6,748         		&    14,270   		     \\ \hline \hline

  \multirow{4}{*}{IndPines} &   \multirow{2}{*}{ \algnameSL} & $\text{acc}_{\text{val}}$     	&	90.2\%  	&      90.6\% 	&    91.2\%    &   91.4\%     \\ \cline{3-7} 
&   	& $\pazocal{C}_{\text{ML}} $       										&	 2,559 	&     2,892  	&    3,718    &   3,967      \\ \cline{2-7}   
&  \multirow{2}{*}{ {\algname }}  & $\text{acc}_{\text{val}}$    	&    90.1\%		&   91.7\%      	&      91.8\%      		&       	91.9\%	     \\ \cline{3-7}   
&   & $\pazocal{C}_{\text{ML}} $    	&    1,662		&      1,220   	&      984      		&     863  		     \\ \hline
\end{tabular}
\caption{The total computational cost $\pazocal{C}_{\text{ML}}$ of \algname \  and \algnameSL \ with respect to increasing number of levels. 
For \textit{IndPines}, we also report the highest achieved $\text{acc}_{\text{val}}$.}
\label{tab:ML}
\end{table}

\subsubsection{\algname \ with and without active-set approach}
As mentioned earlier, the conditions imposed on $R_{\ell}$ and $P_{\ell}$ are extremely general. 
This flexibility allows us to propose yet another variant of the \algname\ algorithm. 
Motivated by the results reported in~\cite{kornhuber1994monotone} and~\cite{kopanivcakova2023multilevel} for linear multigrid and multilevel trust-region methods, respectively, we now consider a variant of \algname\ that incorporates an active set strategy.

In particular, before \algname\ descends to a lower level and performs a recursive iteration (Step 3), we identify the active set using the current iterate $x_{\ell, k}$, i.e., 
\begin{align}
\pazocal{A}_{\ell, k}(x_{\ell, k}) := \{ j \in \{1, \ldots, n_{\ell}  \}  \ | \ l_{\ell, k, j} = x_{\ell, k, j} \ \  \text{or} \ \ u_{\ell, k, j} = x_{\ell, k, j}  \}.
\end{align}
The components of the active set~$\pazocal{A}_{\ell, k}$ are then held fixed and cannot be altered by the lower levels. 
To this end, we define the truncated prolongation operator $\widetilde{P}_{\ell-1}$ as
\begin{align}
(\widetilde{P}_{\ell-1})_{pt} = 
\begin{cases}
0,  \quad &\text{if }\, p \in \pazocal{A}_{\ell, k}(x_{\ell, k}), \\
({P}_{\ell-1})_{pt}, \quad& \text{otherwise}.
\end{cases}
\label{eq:truncation}
\end{align}
Thus, the operator~$\widetilde{P}_{\ell-1}$ is obtained from the prolongation operator~${P}_{\ell-1}$ by simply setting the~$p$-th row of~${P}_{\ell-1}$ to zero for all~$p \in \pazocal{A}_{\ell, k}$. 
As before, we set~$\widetilde{R}_{\ell-1} =  2^{-d} \widetilde{P}^T_{\ell-1}$. 
Note that $(\widetilde{P}_{\ell-1} s_{\ell-1})_p = 0$ for all $p \in \pazocal{A}_{\ell, k}$ and some coarse-level correction~$s_{\ell-1, k} = x_{\ell-1, k} - x_{\ell, 0}$.

These iteration-dependent truncated transfer operators~$\widetilde{P}_{\ell-1}$ and $\widetilde{R}_{\ell-1}$ are then used in place of $P_{\ell-1}$ and $R_{\ell-1}$ during the execution of the recursive iteration (Step 3) of the \algname\ algorithm. 
In addition, we use them to define the coarse-level model\footnote{Non-quadratic coarse-level models can also be used. However, their construction requires the assembly using truncated basis functions, which is tedious in practise~\cite{kovcvara2016first}.} $h_{\ell-1}$, 
by simply restricting the quadratic model from level~$\ell$ to level~$\ell-1$, as follows:
\begin{align}
h_{\ell-1}(x_{\ell-1})  =  (\widetilde{R}_{\ell-1} g_{\ell, k})^T s_{\ell-1} + \frac{1}{2} s_{\ell-1}^T  (\widetilde{R}_{\ell-1} B_{\ell, k} \widetilde{P}_{\ell-1}) s_{\ell-1}.
\label{eq:galerkin_CM}
\end{align}
The model~$h_{\ell-1}$ is used in place of~$f_{\ell-1}$ on the coarse level. 
This construction ensures that the components of $g_{\ell, k}$ and $B_{\ell, k}$ associated with the active set~$\pazocal{A}_{\ell, k}$ are excluded from the definition~$h_{\ell-1}$.

\begin{figure}[t]
\centering
\includegraphics{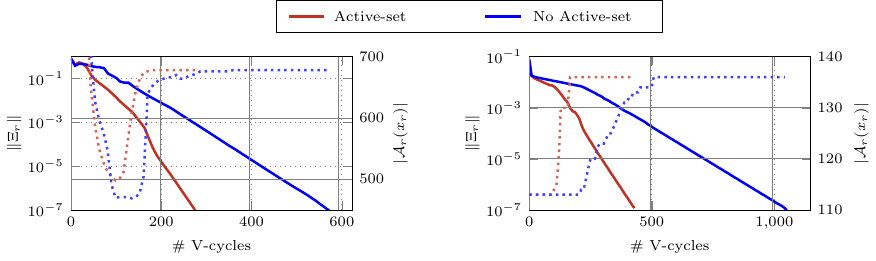}
  \caption{The comparison of the convergence of the \algname \ with and without the active-set strategy. Both variants use the Galerkin-based coarse-level models~\eqref{eq:galerkin_CM}. 
Left: \textit{MinSurf} (three levels). 
Right: \textit{Membrane} (five levels). }
  \label{fig:models_AS_vs_noAS}
\end{figure}

Finally, we also point out that, since~$\widetilde{P}_{\ell-1}$ has fewer nonzero columns than~$P_{\ell-1}$, the resulting lower-level bounds obtained via~\eqref{lower-level-bounds} are less restrictive than those derived using $P_{\ell-1}$. 
This, in turn, allows for larger steps on the coarse level and leads to improved convergence.
To demonstrate this phenomenon, we depict the convergence of the \algname \ algorithm with and without the active set strategy for the \textit{MinSurf} and \textit{Membrane} examples. 
Figure~\ref{fig:models_AS_vs_noAS} illustrates the obtained results. 
As we can see, the variant with the active set strategy converges significantly faster\footnote{The convergence of \algname \ without the active set strategy performs comparably for \textit{MinSurf} and \textit{Membrane} examples regardless of whether $\tau$-corrected or restricted quadratic models are used on the coarser levels.} - approximately by a factor of two.
We also observe that convergence is particularly accelerated once the active set on the finest level~$\pazocal{A}_r$ is correctly identified by the \algname \ algorithm.

\subsubsection{Sensitivity of \algname \ to the choice of the step-size}

\begin{figure}[t]
\centering
\includegraphics{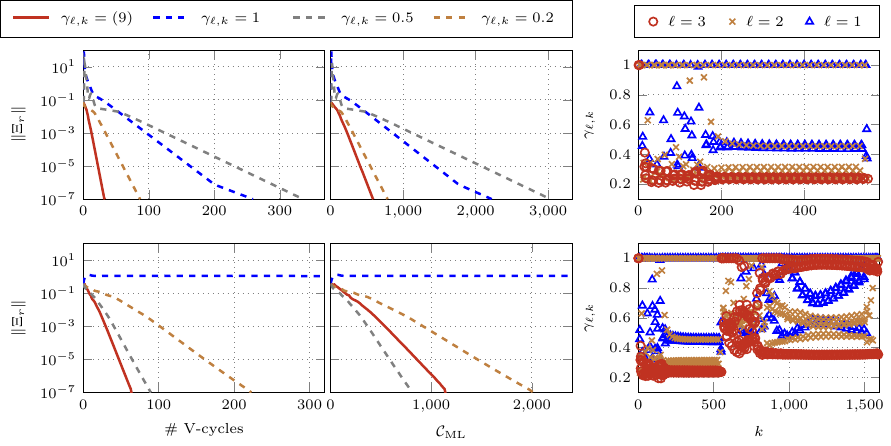}
  \caption{Convergence behavior of \algname \ for \textit{Membrane} (Top) and \textit{MinSurf} (Bottom) examples with three levels. 
Left and Middle: Convergence as a function of V-cycles and computational cost for different choices of~$\gamma_{\ell, k}$. 
The values of~$\gamma_{\ell, k}$, which remain constant across all levels and iterations, are represented by dashed lines. 
Right: The values of~$\gamma_{\ell, k}$ obtained by means of~\eqref{gamma-def} on different levels during the iteration process. 
Here, the symbol $k$ denotes the global iteration counter, incremented across all levels and iteration steps.}
  \label{fig:influence_of_gammas}
\end{figure}

This subsection aims to highlight the importance of carefully selecting the step size~$\gamma_{\ell, k}$ across all levels and iterations. 
To this end, we consider the \textit{Membrane} and \textit{MinSurf} problems and evaluate the performance of the three-level \algname\ algorithm using either an iteration-dependent~$\gamma_{\ell, k}$, computed as in~\eqref{gamma-def}, or constant~$\gamma_{\ell, k}$, choosen from the set~$\{0.2, 0.5, 1.0\}$.
As shown in Figure~\ref{fig:influence_of_gammas} (left), using the iteration-dependent~$\gamma_{\ell, k}$ leads to the fastest convergence in terms of the number of required V-cycles. 
Notably, if fixed~$\gamma_{\ell, k}$ is chosen too large, the algorithm may fail to converge due to violation of the sufficient decrease condition, specified in~\eqref{inside-and-decr}.

Figure~\ref{fig:influence_of_gammas} (right) shows the values of~$\gamma_{\ell, k}$ on different levels, as determined by~\eqref{gamma-def}. 
As we can see, the values of~$\gamma_{\ell, k}$ vary significantly throughout the iteration process and, typically, larger~$\gamma_{\ell, k}$ are determined on lower levels. 
We also note that when accounting for the cost of gradient evaluation, it is possible to identify a constant~$\gamma_{\ell, k}$ for which the overall computational cost~$\pazocal{C}_{\text{ML}}$  is comparable to that achieved using iteration-dependent step sizes. 
This is due to the fact that evaluating~$(s_{\ell,k}^L)^T B_{\ell,k} s_{\ell,k}^L$ via complex-step finite differences requires an additional gradient evaluation on each iteration.  
We emphasize, however, that this comparable~$\pazocal{C}_{\text{ML}}$ comes at the expense of costly hyper-parameter tuning required to determine an effective~$\gamma_{\ell, k}$. 
Therefore, we restrict the use of fixed step sizes only to the DDN examples, as is common practice in machine learning.

\subsubsection{Sensitivity of \algname\  to noise}
\begin{figure}[t]
\centering
\includegraphics{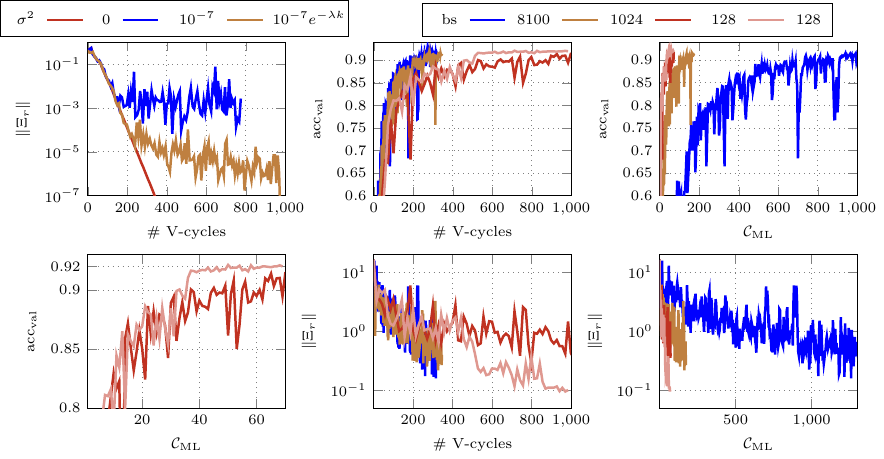}
\caption{\textit{Top Left:} Convergence of the \algname \ without noise (red), with constant noise (blue, $\sigma^2 =10^{-7}$) and the noise reduced using the exponential scheduler (brown) for \textit{MinSurf} (two levels). 
\textit{Bottom Left:} Convergence of the \algname \ with noise due to gradients subsampled with batch size (bs) of $128$ for \textit{IndPines} example (three levels). 
The red color represents results with constant noise, while pink color denotes results with learning rate (lr) scheduler, where lr is dropped by factor of $10$ at $\pazocal{C}_{\text{ML}}$ equal to $35$ and $55$. 
\textit{Middle/Right:} Convergence of the \algname \ with noise due to gradients subsampled for \textit{IndPines} example. 
The blue color stands for exact gradient evaluations, i.e., with the whole dataset (bs=8,100), while brown (bs=1,024) and red/pink (bs=128) consider inexact gradient evaluations.}
\label{fig:noise_ML}
\end{figure}

The main motivation for developing and using OFFO algorithm lies in their performance in the presence of noise.
To illustrate the impact of noise on the convergence behaviour of \algname, we first consider the \textit{MinSurf} problem. 
Figure~\ref{fig:noise_ML} (top left) displays the convergence of the algorithm in the absence of noise (red line) and under constant additive noise (blue line). 
Thus, $g_{\ell, k} = \nabla f(x_{\ell, k}) + \epsilon_{\ell, k}$, where each component of random noise vector $\epsilon_{\ell, k}$ is obtained from $\pazocal{N}(0, \sigma_{\ell, k}^2)$, with $\sigma_{\ell, k}^2$ denoting the noise variance.
Here, we set ~$\sigma_{\ell, k}^2 = 10^{-7}$, for all $(\ell, k)$.
As observed, the presence of noise does not significantly affect the convergence until $\|\Xi_{r}\|$ becomes comparable to the noise level. 
This is particularly beneficial for applications that do not require highly accurate solutions, such as the training of DNNs.
However, if a high-accuracy solution is required, it can be obtained by progressively reducing the noise in the gradient evaluations. 
Following standard practice, one may employ a noise reduction scheduler. 
For example, we can employ exponential decay scheme~\cite{bottou2018optimization}, with $\sigma_{\ell, k}^{2} = \sigma_{\ell, 0}^{2}  e^{-\lambda k}$, where we set~$\sigma_{\ell, 0}^{2}  = 10^{-7}$ and $\lambda = 5 \cdot 10^{-2}$.
As illustrated by the brown line in Figure~\ref{fig:noise_ML} (top left), incrementally reducing the noise enables \algname\  to converge to a critical point with user prescribed tolerance.
Note that the choice of noise reduction scheduler can have a significant impact on the convergence speed.

Understanding \algname's   sensitivity to noise is particularly relevant in the case of DNNs, where inexact gradients, obtained by subsampling, are often used to reduce the iteration cost.
Figure~\ref{fig:noise_ML} (middle and right) demonstrates the convergence of \algname\ for the {\textit{IndPines}} example, when different numbers of samples are used to evaluate the gradient.
In particular, we set the batch size ({bs}) to 8,100, 1,024, and 128, where 8,100 corresponds to the full dataset.
As we can see from the obtained results, \algname \ with the full dataset converges the fastest in terms of V-cycles.
However, since the cost of evaluating the gradient scales linearly with the number of samples in the batch, \algname \  with {bs}~$=128$ is the most cost-efficient.
Looking at $\text{acc}_{\text{val}}$ and $\|\Xi_{r}\|$, we observe that \algname\ with bs=$8,100$ produces the most accurate DNN.
As before, the noise introduced by subsampling can be mitigated during the training process. 
To this end,  we use a learning-rate step decay scheduler\footnote{As an alternative to reducing the learning rate, one may gradually increase the mini-batch size to mitigate the noise in gradient evaluations introduced by subsampling.}~\cite{goodfellow2016deep}, in which the learning rate is reduced by a factor of ten at $\pazocal{C}_{\text{ML}}$ equal to $35$ and $55$.
As shown in Figure~\ref{fig:noise_ML} (bottom left), this strategy enhances the quality of the resulting model while maintaining the low computational cost per iteration associated with small-batch gradient evaluations.
Moreover, in the long run, the final model accuracy surpasses that achieved by using the whole dataset, highlighting the role of noisy gradients in the early stages of training for improved generalization.

\subsection{Numerical performance of \algnamed \ }
\label{sec:DD_numerics}

In this section, we evaluate the numerical performance of \algnamed.
To this aim, we estimate its total parallel computational cost as
\begin{align}
\pazocal{C}_{\text{DD}} = \bigg(\sharp_r \pazocal{C}_r \bigg) + \bigg(\frac{n^{(p)}}{n_{r}} \sharp^{(p)} \pazocal{C}_0 \bigg),
\label{eq:cost_DD}
\end{align}
where~$n^{(p)}$ denotes the number of degrees of freedom in the largest subdomain, and $\sharp^{(p)}$ represents the number of steps performed on all subdomains\footnote{In our tests, the same number of Taylor iterations is performed in all subdomains.}.
Thus, the first term in~\eqref{eq:cost_DD} accounts for the cost of the Taylor iteration, while the second term corresponds to the cost of subdomain iterations, which can be carried out in parallel.

Table~\ref{tab:DD_Taylor} reports the results obtained for the \textit{Membrane} and \textit{MinSurf} examples with respect to increasing number of subdomains and overlap.
As we can observe, the parallel computational cost~$\pazocal{C}_{\text{DD}}$ decreases with increased number of subdomains, reflecting the benefits of possible parallelization.
Furthermore, increasing the overlap size contributes to an additional reduction in the computational cost.
It is worth noting that the observed parallel speedup of \algnamed \ compared to \adagrad \  is below the ideal.
For example, for eight subdomains and an overlap size of two, the speedup is approximately of a factor of four.
We attribute this suboptimal theoretical scaling primarily to the cost associated with serial Taylor iterations.
As we will discuss in Section~\ref{sec:hyb_algo}, this limitation can be addressed by developing hybrid multilevel-domain-decomposition algorithms, which effectively combine \algname\ and \algnamed.

\begin{figure}[t]
\centering
\includegraphics{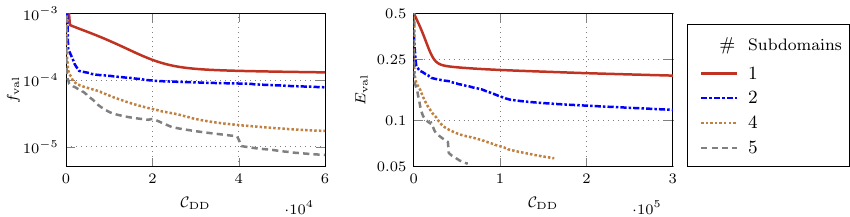}
  \caption{The parallel computational cost~$\pazocal{C}_{\text{DD}}$ of \adagrad ($\#$ subdomains=1) and \algnamed \ for \textit{Aniso} with varying number of subdomains (no overlap).
  Left: The validation loss $f_{\text{val}}$. Right: The relative validation error ${E_{\text{val}} = 1/n_s \sum_{s=1}^{n_s} (y_s^D - y_s) / \| y_s \| }$, where $y_s^D$ is the DeepONet-inferred solution and $y_s$ is the target solution.
  The gradient evaluations are performed using batch-size $bs=1,062$.}
  \label{fig:DD_machine_learning}
\end{figure}

Next, we assess the performance of \algnamed\ for training DeepONet.
Figure~\ref{fig:DD_machine_learning} illustrates the results obtained for the \textit{Aniso} example.
Consistent with the previous observations, \algnamed\ achieves a significant speedup compared to standard \adagrad.
Moreover, as for FEM examples, increasing the number of subdomains enables greater parallelism, thereby reducing the parallel computational cost~$\pazocal{C}_{\text{DD}}$.
Notably, the use of \algnamed\ also yields more accurate DNNs, as evidenced by lower relative errors in the resulting DeepONet predictions.
This improvement in accuracy is particularly important, as the reliability of the DeepONet plays a crucial role in a wide range of practical applications.

\begin{table}[htbp]
\centering
\begin{tabular}{|l||r|r|r||r|r|r|r|r|r|r|r|}
\hline
\multirow{2}{*}{$\#$ Subdomains/Overlap} & \multicolumn{3}{c||}{Membrane}  &  \multicolumn{3}{c|}{MinSurf}  \\ \cline{2-7}
								 &   \multicolumn{1}{c|}{\textbf{0}} & \multicolumn{1}{c|}{\textbf{2}} & \multicolumn{1}{c||}{\textbf{4}}  & \multicolumn{1}{c|}{\textbf{0}} & \multicolumn{1}{c|}{\textbf{2}} & \multicolumn{1}{c|}{\textbf{4}}  \\ \hline\hline
\textbf{1} (\adagrad)          & 36,054      & --         & --    & 41,136     & --         & --       \\ \hline
\textbf{2}                       	&   21,104    	& 20,340      	&  20,230    	&    22,958     	&    21,906        &  21,868    \\ \hline
\textbf{4}                       	&    12,864   	&  12,198    	&  12,204     	&     14,162       &   13,538         &  12,952    \\ \hline
\textbf{8}                       	&      9,286  	&   8,574    	&   8,490   	&     9,424      	&     9,132      	&  9,030      \\ \hline
\end{tabular}
\caption{The parallel computational cost~$\pazocal{C}_{\text{DD}}$ of \algnamed \ and
\adagrad \  for \textit{Membrane} (four refinement levels) and  \textit{MinSurf} (four refinement levels) examples. }
\label{tab:DD_Taylor}
\end{table}

\subsection{Numerical performance of \algnamehyb \ }
\label{sec:hyb_algo}

As discussed in Section~\ref{sec:DD_numerics}, the main benefit of the \algnamed\ algorithm is that it enables the parallelization.
However, the cost associated with the sequential Taylor step prohibits the ideal speedup in practice.
To improve parallelization capabilities, we draw inspiration from established practices in the domain-decomposition literature~\cite{dolean2015introduction}, where coarse spaces are known to enhance the scalability and robustness of the algorithms at low computational cost.
To this end, we now numerically evaluate the effectiveness of a hybrid algorithm, named \algnamehyb, which effectively combines the \algname\ and \algnamed\ methods.
In particular, we replace the Taylor iteration (Step 2) in Algorithm~\ref{rec-algo-d} with a call to the \algname\ method, i.e., Algorithm~\ref{rec-algo}.
Here, \algname\ is configured with two levels, where the coarse level is constructed by discretizing the original problem on a mesh coarsened by a factor of eight. 
This coarse-level approximation is then adjusted using the tau-correction approach, as defined in~\eqref{tau-corrected}.
Our implementation then alternates between performing ten coarse-level iterations and ten decomposition steps.
Note, if~\eqref{nogain} is satisfied, the recursive steps are replaced with a Taylor iteration on $f_r$. 

The parallel computational cost of our hybrid \algnamehyb \ method is given as 
\begin{align}
\pazocal{C}_{\text{ML-DD}} =  \sum_{\ell=1}^r \frac{n_{\ell}}{n_{r}} \sharp_\ell \, \pazocal{C}_r + \bigg(\frac{n^{(p)}}{n_{r}} \sharp^{(p)} \pazocal{C}_0 \bigg).
\end{align}
Here,  we set~$r=2$, i.e.,~$f_{r}=f_2=f$, $f_1$ is associated with the coarse-level discretization, and $f_0$ is the extended model associated with the subdomains computations, given by~\eqref{f0-def}. 

Table~\ref{tab:ML_DD} reports the results obtained for the \textit{Membrane} and \textit{MinSurf} examples.
Firstly, we observe that even \algnamed\ algorithm exhibits algorithmic scalability, i.e.,~the number of V-cycles remains constant with an increasing number of subdomains.
This behavior can be observed due to our specific implementation, which alternates ten subdomain steps with one global Taylor iteration.
However, as we can see, replacing a Taylor iteration with ten coarse-level iterations leads to a significant improvement in terms of convergence, as indicated by the reduced number of V-cycles.
Moreover, the parallel computational cost of \algnamehyb\ is substantially lower than that of \algnamed, demonstrating strong potential for our envisioned future work, focusing on integrating \algnamehyb\ into parallel finite-element and/or deep-learning frameworks.

\begin{table}[htbp]
\centering
\begin{tabular}{|l||r|r|r|r||r|r|r|r|r|r|r|r|r|r|}
\hline
\multirow{3}{*}{$\#$  Subd.} & \multicolumn{4}{c||}{Membrane}  &  \multicolumn{4}{c|}{MinSurf} \\ \cline{2-9}
                                 			& \multicolumn{2}{c|}{\algnamed}  &  \multicolumn{2}{c||}{\algnamehyb}   & \multicolumn{2}{c|}{\algnamed}  &  \multicolumn{2}{c|}{\algnamehyb}   \\ \cline{2-9} 
                                 			& \multicolumn{1}{c|}{$\pazocal{C}_{\text{DD}}$} &  \multicolumn{1}{c|}{V-cycles}& \multicolumn{1}{c|}{$\pazocal{C}_{\text{ML-DD}}$} &  \multicolumn{1}{c||}{V-cycles}& \multicolumn{1}{c|}{$\pazocal{C}_{\text{DD}}$} &  \multicolumn{1}{c|}{V-cycles}& \multicolumn{1}{c|}{$\pazocal{C}_{\text{ML-DD}}$} &  \multicolumn{1}{c|}{V-cycles} \\ \hline\hline
						
\textbf{1} (\adagrad)      	& 36,054 		& 36,054 		& 36,054   	&36,054 	& 41,136		& 41,136	& 41,136 	& 41,136       \\ \hline 
\textbf{2}                       	&   20,340 	&1,784		&      	 586 		&57		&  21,906  	&1,826   	&  3,864 	&370    \\ \hline
\textbf{4}                       	&    12,198    	&1,900		& 	 312 		&58		&  13,538   	&1,934     &  2,024 	&385  \\ \hline
\textbf{8}                       	&     8,574   	&2,158		& 	 168 		&60		&  9,132   		&2,032     &  1,142 	&406  \\ \hline
\textbf{16}                     	&      6,318  	&1,952		&   	 92 		&60		&   7,806   	&2,403     &  612 	&405  \\ \hline
\end{tabular}
\caption{The parallel computational cost and the number of V-cycles of \algnamed, and  \algnamehyb \  with overlap equal to two for \textit{Membrane} (four refinement levels) and  \textit{MinSurf} (four refinement levels). }
\label{tab:ML_DD}
\end{table}

\section{Conclusions and perspectives}
\label{section:conclusion}

We have proposed an OFFO algorithmic framework inspired by AdaGrad which, unlike the original 
formulation, can handle bound constraints on the problem variables and incorporate 
curvature information when available. Furthermore, the framework is designed to 
efficiently exploit the structural properties of the problem, including hierarchical multilevel 
decompositions and standard additive Schwarz domain decomposition strategies.
The convergence of this algorithm has then been studied in the stochastic setting allowing 
for noisy and possibly biased gradients, and its evaluation complexity for computing 
an $\epsilon$-approximate first-order critical point was proved to be 
$\calO(\epsilon^{-2})$ (without any logarithmic term) with high probability.

Extensive numerical experiments have been discussed, showing the significant computational 
gains obtainable when problem structure is used, both in the hierarchical and domain 
decomposition cases, but also by combining the two approaches.
These gains are demonstrated on problems arising from discretized PDEs and machine 
learning applications such as the training of ResNets and DeepONets. The use of cheap 
but accurate gradient-based finite-difference approximations for curvature information 
was also shown to be crucial for numerical efficiency on the PDE-based problems.

While already of intrinsic interest, these results are also viewed by the authors as a step towards 
efficient, structure exploiting methods for problems involving more general constraints. Other  
theoretical and practical research perspectives include, for instance, the use of momentum, 
active constraints' identification, adaptive noise reduction scheduling and the 
streamlining of strategies mixing domain decomposition and hierarchical approaches, 
in particular in the domain of deep neural networks.

{\footnotesize
\section*{\footnotesize Acknowledgements}
This work benefited from the AI Interdisciplinary Institute ANITI, funded by the 
France 2030 program under Grant Agreement No. ANR-23-IACL-0002. Moreover, the numerical 
results were carried out using HPC resources from GENCI-IDRIS (Grant No. AD011015766).
Serge Gratton and Philippe Toint are grateful to Defeng Sun, Xiaojun Chen and Zaikun Zhang 
of the Department of Applied Mathematics of Hong-Kong Polytechnic University for their 
support during a research visit in the fall 2024.

\bibliography{./gkt3.bib}
\bibliographystyle{plain}
} 

\end{document}